\documentclass[preprint, 10pt]{elsarticle}

\usepackage[margin=2.5cm]{geometry}
\usepackage{setspace}
\usepackage{lmodern}
\usepackage{amsmath,amssymb}
\usepackage{amsfonts,amsthm}
\usepackage{lineno}
\usepackage{graphicx}
\usepackage{bm, bbm}
\usepackage{algorithm}
\usepackage{physics}
\usepackage{color}
\usepackage{pxfonts}
\usepackage[footnotesize,bf]{caption}
\usepackage{subcaption}
\usepackage{mathtools}
\usepackage{hhline}
\usepackage[T1]{fontenc}
\usepackage{multirow}
\usepackage{overpic}
\usepackage{hyperref}
\usepackage{placeins}

\date{}

\newtheorem{remark}{Remark}[section]

\graphicspath{{figure/}}
\begin{document}
	\begin{frontmatter}
		\title{Latent representation learning based model correction and uncertainty quantification for PDEs}
		\author[SHNU]{Wenwen Zhou}
		
		\author[UIC]{Xiaodong Feng}
		\author[SHNU]{Ling Guo}
		\author[SJTU]{Hao Wu} 
	
		\vspace{0.2cm}
		\address[SHNU]{Department of Mathematics, Shanghai Normal University, Shanghai, China.}
		
		\vspace{0.2cm}
		\address[UIC]{Faculty of Science and Technology, Beijing Normal-Hong Kong Baptist University, Zhuhai 519087, China.}
		\vspace{0.2cm}
		\address[SJTU]{School of Mathematical Sciences, Institute of Natural Sciences, and MOE-LSC, Shanghai Jiaotong University, Shanghai, China.}
		
		\begin{abstract}
			Model correction is essential for reliable PDE learning when the governing physics is misspecified due to simplified assumptions or limited observations. In the machine learning literature, existing correction methods typically operate in parameter space, where uncertainty is often quantified via sampling or ensemble-based methods, which can be prohibitive and motivates more efficient representation-level alternatives.			
			To this end, we develop a latent-space model-correction framework by extending our previously proposed LVM-GP solver, which couples latent-variable model with Gaussian processes (GPs) for uncertainty-aware PDE learning. Our architecture employs a shared confidence-aware encoder and two probabilistic decoders, with the solution decoder predicting the solution distribution and the correction decoder inferring a discrepancy term to compensate for model-form errors.
			The encoder constructs a stochastic latent representation by balancing deterministic features with a GP prior through a learnable confidence function. Conditioned on this shared latent representation, the two decoders jointly quantify uncertainty in both the solution and the correction under soft physics 
			constraints with noisy data. An auxiliary latent-space regularization is introduced to control the learned representation and enhance robustness. This design enables joint uncertainty quantification of both the solution and the correction within a single training procedure, without parameter sampling or repeated retraining. Numerical experiments show accuracy comparable to Ensemble PINNs and B-PINNs, with improved computational efficiency and robustness to misspecified physics.
		\end{abstract}
	\end{frontmatter}
	
	\section{Introduction}\label{sec:introduction}
	
	Accurate partial differential equation (PDE) models are central to scientific computing and engineering, providing the foundation for simulation-based prediction, analysis, and control of complex physical systems \cite{lapidus1999numerical, leveque2007finite, karniadakis2021physics}. In practice, the governing physics is often simplified \cite{oberkampf2002error, box1979robustness, benner2017model} and the available data are limited \cite{raissi2019physics}. As a result, model form errors and operator misspecification are difficult to avoid \cite{subramanian2021model}. These issues are particularly pronounced for nonlinear dynamics \cite{thompson2002nonlinear}: small structural inaccuracies can propagate and amplify over long time horizons \cite{hagstrom2009complete, calvo2011error} or in high-dimensional regimes \cite{barbier2019optimal, krishnapriyan2021characterizing}, resulting in substantial discrepancies between numerical predictions and real-world observations \cite{lorenz1989computational, oberkampf2006measures}. Model correction has therefore become an important research topic for PDE-based modeling \cite{kennedy2001bayesian, zou2024correcting}, especially when solution fidelity is critical for practical applications \cite{oberkampf2002error, oberkampf2010verification}.
	
	Classical model correction typically relies on parametric calibration \cite{ tarantola2005inverse}, closure modeling \cite{pope2001turbulent, wilcox1998turbulence, pan2018data}, residual or bias correction \cite{kennedy2001bayesian}, where a low-dimensional discrepancy model is tuned to reconcile simulations with measurements. Such strategies can be effective when dominant error sources \cite{oberkampf2002error} are well understood and the discrepancy admits a compact parameterization \cite{sargsyan2019embedded}. However, they can be restrictive when misspecification is structural \cite{brynjarsdottir2014learning,morrison2018representing}, spatially varying \cite{xu2015bayesian}, and when observations provide only partial coverage of the state space \cite{evensen2009data}.
	
	Recent machine-learning-based approaches increase flexibility by representing model discrepancies with data-driven surrogate models \cite{raissi2019physics, brunton2020machine}. A widely used direction is additive discrepancy correction \cite{vagnoli2025local, duraisamy2019turbulence}, which augments the governing equations with a learned discrepancy term to compensate for systematic bias; representative examples include hybrid forecasting for bias reduction in weather prediction \cite{farchi2021using} and PINN-framework correction layers that learn discrepancy terms jointly with the solution \cite{zou2024correcting}. Beyond additive discrepancy correction, increasing attention has been devoted to differential operator correction, where one modifies the differential operator itself to account for operator level misspecification, including symbolic and sparse regression formulations \cite{brunton2016discovering, zhang2024discovering, podina2022pinn, cranmer2023interpretable} for discovering corrected operators. While these methods can substantially improve predictive fidelity, the correction task is often weakly identifiable \cite{kennedy2001bayesian, tarantola2005inverse, kaipio2005statistical} under noisy or sparse measurements: distinct corrected models may fit the data comparably well, yielding ambiguity in both the inferred correction and the resulting PDE predictions. This ambiguity makes uncertainty quantification (UQ) \cite{smith2024uncertainty, psaros2023uncertainty} an essential component of model correction. Existing UQ treatments based on ensembles or Bayesian inference in parameter space can be computationally expensive \cite{psaros2023uncertainty, yang2021b, pensoneault2024efficient} and may struggle to represent uncertainty induced by differential operator misspecification \cite{zou2024neuraluq, shen2023picprop}.
	
	In this work, we develop an uncertainty-aware PDE model correction framework built on the latent variable model-Gaussian process (LVM-GP) \cite{feng2025lvm}. The key design is a shared confidence-aware latent representation that simultaneously drives the solution model and the correction model, coupling their epistemic uncertainty through the same stochastic latent field. The solution and the correction are learned jointly using a unified objective that combines data likelihood terms with physics-based constraints \cite{raissi2019physics}, together with a latent regularization that prevents degeneracy of the  encoder in latent space \cite{guo2024ib, kingma2013auto}. This construction yields coherent uncertainty propagation across the solution, the corrected forcing, and the inferred discrepancy, while retaining the useful inductive bias of the possibly misspecified physics. Our main contributions are summarized as follows:
	
	\begin{itemize}
		\item We introduce a latent-space model correction framework in which the solution and the correction share a confidence-aware stochastic latent field. This shared representation enables consistent uncertainty propagation across all quantities of interest.
		
		\item We develop a joint learning strategy that embeds PDE information directly into the data likelihood. As a result, the solution and the correction are learned simultaneously from a single coupled objective. We further introduce a latent regularization term to stabilize training and prevent collapse of the stochastic latent component.
		
		\item We validate the proposed framework on several benchmark problems, including an ODE system, a reaction-diffusion equation, and two-dimensional channel and cavity flows. The numerical results show improved accuracy and uncertainty calibration in the presence of noisy data and model misspecification.
	\end{itemize}
	
	The remainder of this paper is organized as follows. Section~\ref{sec:preliminaries} reviews model correction ideas in physics-informed neural networks (PINNs) under misspecified operators and summarizes the LVM-GP framework for uncertainty-aware PDE learning. Section~\ref{Methodology} presents the proposed LVM-GP model correction method, detailing the shared latent representation, the neural-operator decoder, and the latent-conditioned correction mechanism, together with the associated training and uncertainty quantification procedure. Section~\ref{Experiments and Results} demonstrates the performance of the proposed method through several numerical examples. Section~\ref{Summary} concludes the paper and discusses limitations and future directions.
	
	\section{Preliminaries}\label{sec:preliminaries}
	
	In this section, we briefly review recent approaches to model correction in PINNs for handling misspecified PDE operators. We then introduce uncertainty-aware PDE solvers based on the LVM-GP model.

	\subsection{Model corrections in physics-informed neural networks}\label{subsec:pde-data-correction}
	
	We consider the following PDE defined on a domain $\Omega \subset \mathbb{R}^{d}$ with boundary $\partial\Omega$:
	\begin{equation}\label{eq:general_PDE}
		\begin{aligned}
			\mathcal{N}_{\bm{x}}[u](\bm{x}) &= f(\bm{x}), \qquad \bm{x}\in\Omega,\\
			\mathcal{B}_{\bm{x}}[u](\bm{x}) &= b(\bm{x}), \qquad \bm{x}\in\partial\Omega,
		\end{aligned}
	\end{equation}
	where $u$ is the unknown solution, $\mathcal{N}_{\bm{x}}$ is a differential operator, and $\mathcal{B}_{\bm{x}}$ is the boundary operator.
	
	We assume that the available data consist of noisy pointwise measurements of $u$, $f$, and $b$:
	\begin{equation*}
		\mathcal{D}=\mathcal{D}_{u}\cup \mathcal{D}_{f}\cup \mathcal{D}_{b},
	\end{equation*}
	with
	\begin{equation*}
		\begin{aligned}
			\bar{u}^{(i)} &= u(\bm{x}^{(i)}_{u}) + \epsilon^{(i)}_{u}, \quad i = 1,\dots,N_{u},\\
			\bar{f}^{(i)} &= f(\bm{x}^{(i)}_{f}) + \epsilon^{(i)}_{f}, \quad i = 1,\dots,N_{f},\\
			\bar{b}^{(i)} &= b(\bm{x}^{(i)}_{b}) + \epsilon^{(i)}_{b}, \quad i = 1,\dots,N_{b},
		\end{aligned}
	\end{equation*}
	where $\epsilon^{(i)}_{u},\epsilon^{(i)}_{f},\epsilon^{(i)}_{b}$ are independent zero-mean Gaussian noises.
	
	The main idea of PINNs is to approximate the unknown solution field $u(\bm{x})$ by a parametrized neural network
	$u_{\mathrm{NN}}(\bm{x};\theta)$. Given the data sets $\mathcal{D}_u$, $\mathcal{D}_f$, and $\mathcal{D}_b$, a standard PINN is trained by minimizing the following least-squares loss function:
	\begin{equation*}\label{eq:pinn_loss}
		\mathcal{L}(\theta)= \omega_f\,\mathcal{L}_f(\theta) + \omega_b\,\mathcal{L}_b(\theta) + \omega_u\,\mathcal{L}_u(\theta),
	\end{equation*}
	where
	\begin{equation*}\label{eq:pinn_loss_terms}
		\mathcal{L}_f(\theta)=\frac{1}{2}\sum_{i=1}^{N_f}\Big|
		\mathcal{N}_{\bm{x}}[u_{\mathrm{NN}}](\bm{x}_f^{(i)})-\bar f^{(i)}
		\Big|^2,\qquad
		\mathcal{L}_b(\theta)=\frac{1}{2}\sum_{i=1}^{N_b}\Big|
		\mathcal{B}_{\bm{x}}[u_{\mathrm{NN}}](\bm{x}_b^{(i)})-\bar b^{(i)}
		\Big|^2,
	\end{equation*}
	and
	\begin{equation*}\label{eq:pinn_data_term}
		\mathcal{L}_u(\theta)=\frac12\sum_{i=1}^{N_u}\Big|
		u_{\mathrm{NN}}(\bm{x}_u^{(i)};\theta)-\bar u^{(i)}
		\Big|^2.
	\end{equation*}
	Here $\{\omega_f,\omega_b,\omega_u\}$ are prescribed weights balancing the contributions of the PDE residual, boundary condition and data mismatch terms.
	
	When the governing operator is misspecified, the resulting optimization problem is inherently inconsistent: the data-misfit and residual terms generally cannot be minimized simultaneously. In particular, the residual-based loss is dominated by model discrepancy, which biases the learned solution toward satisfying the misspecified PDE instead of matching the true physics implied by the data. This mismatch leads to a trade-off and can cause biased predictions and miscalibrated uncertainty.
	
	A standard remedy is to augment the misspecified operator by a discrepancy term. Let $\widetilde{\mathcal{N}}_{\bm{x}}$ denote a given misspecified differential operator. Following \cite{zou2024correcting}, one can define the model discrepancy
	\[
	s(\bm{x}) = f(\bm{x}) - \widetilde{\mathcal{N}}_{\bm{x}}[u](\bm{x}), \quad \bm{x}\in\Omega.
	\]
	In practice, $s(\bm{x})$ can be represented by another  neural network and is learned together with $u_{\mathrm{NN}}$ utilizing the same loss structure based on misfit data and residuals. This formulation separates the \emph{known} part of the operator from an \emph{unknown} correction, and provides a convenient way to account for model-form error within the PINN setting.
	
	The same idea can be stated at the operator level. Rather than viewing $s(\bm{x})$ purely as an additive correction to the forcing term, one may introduce a correction operator $\mathcal{C}_{\bm{x}}$ and write
	\begin{equation}\label{eq:operator_correction}
		\big(\widetilde{\mathcal{N}}_{\bm{x}} + \mathcal{C}_{\bm{x}}\big)[u](\bm{x}) = f(\bm{x}), \qquad \bm{x}\in\Omega,
	\end{equation}
	where $\widetilde{\mathcal{N}}_{\bm{x}}$ denotes the misspecified operator and $\mathcal{C}_{\bm{x}}$ is to be identified from data and physiccal constraints; see, e.g., \cite{cao2023residual}. This operator level viewpoint is natural in PDE modeling, since it directly targets the relation between $u(\bm{x})$ and its derivatives encoded by the differential operator.
	
	In this setting, learning $\mathcal{C}_{\bm{x}}$ (or $s(\bm{x})$) from finite and noisy data is typically ill-posed: the correction may not be uniquely identifiable, and small perturbations in the data can lead to noticeable changes in the recovered correction. It is therefore crucial to quantify uncertainty in both the correction and the induced solution, rather than reporting a single point estimate.
	Moreover, many existing uncertainty-aware correction pipelines parameterize correction uncertainty and solution uncertainty separately, which can hinder consistent uncertainty propagation and complicate interpretation and calibration. To address these challenges, we introduce LVM-GP, an uncertainty-aware PDE solver that couples a latent variable model with a Gaussian process prior to produce predictive distributions with calibrated uncertainty.
	\subsection{Latent Variable Model-Gaussian Process (LVM-GP): Uncertainty-Aware PDE Solver}
	LVM-GP is a probabilistic, physics-informed framework for PDE solving with uncertainty quantification \cite{feng2025lvm}. It fuses noisy measurements with the governing PDE to capture both aleatoric and epistemic uncertainties, yielding calibrated predictive distributions for the solution field. The framework consists of two main components: a confidence-aware encoder and a probabilistic decoder.
	
	\begin{itemize}
		\item \textbf{Confidence-Aware Encoder}. The confidence-aware encoder constructs a latent representation by interpolating between a deterministic feature map and a Gaussian-process (GP) prior. Specifically,
		\[
		\bm{z}_{1}(\bm{x},\bm{\omega})
		= \mathrm{diag}(\bm{m}(\bm{x}))\,\bar{\bm{z}}(\bm{x})
		+ \mathrm{diag}(1-\bm{m}(\bm{x}))\,\bm{z}_{0}(\bm{x},\bm{\omega}),
		\]
		where $\bm{m}(\bm{x})\in[0,1]^M$ and $\bar{\bm{z}}(\bm{x})\in\mathbb{R}^M$ are outputs of the encoder network. Here, $\bm{m}(\bm{x})$ serves as a confidence function, while $\bar{\bm{z}}(\bm{x})$ represents the deterministic feature. The stochastic term $\bm{z}_{0}(\bm{x},\bm{\omega})$ is induced by a vector-valued GP prior
		\[
		\bm{z}_{0} \sim \mathcal{GP}(0,K(\bm{x},\bm{x}')),
		\]
		with matrix-valued kernel $K(\bm{x},\bm{x}')\in\mathbb{R}^{M\times M}$. The confidence function $\bm{m}(\bm{x})$ controls the interpolation between the deterministic feature $\bar{\bm{z}}(\bm{x})$ and the GP prior, thereby balancing deterministic and stochastic information in the latent space.
		\item \textbf{Probabilistic Decoder}. Conditioned on the latent variable $\bm{z}_{1}(\bm{x},\bm{\omega})$, a neural operator decoder predicts the solution mean $\mu_u(\bm{x})$. The model further represents the solution field $u(\bm{x})$ and the physical terms $f(\bm{x})$ and $b(\bm{x})$ as conditionally Gaussian random fields, enabling joint uncertainty quantification:
		\[
		\begin{aligned}
			u(\bm{x}) &\sim \mathcal{N}\big(\mu_{u}(\bm{x}), \sigma^{2}_{u}(\bm{x})\big), \\
			f(\bm{x}) &\sim \mathcal{N}\big(\mathcal{N}_{\bm{x}}[\mu_{u}](\bm{x}), \sigma^{2}_{f}(\bm{x})\big), \\
			b(\bm{x}) &\sim \mathcal{N}\big(\mathcal{B}_{\bm{x}}[\mu_{u}](\bm{x}), \sigma^{2}_{b}(\bm{x})\big),
		\end{aligned}
		\]
		where $\mathcal{N}(\cdot,\cdot)$ denotes the Gaussian distribution, and $\sigma^{2}_{u}(\bm{x})$, $\sigma^{2}_{f}(\bm{x})$, and $\sigma^{2}_{b}(\bm{x})$ are variance functions that can be prescribed or parameterized and learned from the data. These variances capture the aleatoric uncertainty associated with observations of $u(\bm{x})$, $f(\bm{x})$, and $b(\bm{x})$, respectively.
	\end{itemize}
	The model parameters can be estimated by maximizing the likelihood of the observed data while regularizing the latent space. Equivalently, we minimize
	\begin{equation}\label{eq:total loss}
		\mathcal{L} = \mathcal{L}_{\mathrm{data}} - \beta\,\mathcal{L}_{\mathrm{reg}} .
	\end{equation}
	Here, $\mathcal{L}_{\mathrm{data}}$ is the negative log-likelihood associated with the predictive distributions of $u$, $f$, and $b$,
	\begin{equation*}
		\mathcal{L}_{\mathrm{data}}
		= \omega_{u}\,\mathcal{L}_{\mathrm{data},u}
		+ \omega_{f}\,\mathcal{L}_{\mathrm{data},f}
		+ \omega_{b}\,\mathcal{L}_{\mathrm{data},b},
	\end{equation*}
	with prescribed weights $\omega_u$, $\omega_f$, and $\omega_b$.
	To prevent the latent variable from becoming unstructured and to keep it close to a prescribed prior, we employ a KL regularizer,
	\begin{equation}\label{eq:reg loss}
		\mathcal{L}_{\mathrm{reg}}
		= \mathbb{E}_{\bm{x} \sim \mathcal{U}(\Omega)}
		\big[ D_{\mathrm{KL}}\big(q_E(\bm{z}\mid \bm{x})\,\|\,p(\bm{z})\big) \big],
	\end{equation}
	where the prior is $p(\bm{z})=\mathcal{N}(0,\mathbf{I}_{M})$ by default. For the confidence-aware encoder, the induced distribution takes the form
	\[
	q_E(\bm{z}\mid \bm{x})=\mathcal{N}\big(\mathrm{diag}(\bm{m}(\bm{x}))\,\bar {\bm{z}}(\bm{x}),\ \mathrm{diag}(1-\bm{m}(\bm{x}))^{2}\big).
	\]
	
	Once trained, LVM-GP generates solution samples by drawing latent variables from the encoder and propagating them through the decoder. The predicted solution is given by the sample mean, and the uncertainty is quantified by the sample variance.
	
	Despite its practical advantages in producing uncertainty-aware predictions, the LVM-GP formulation implicitly assumes that the governing PDE is correctly specified. In many scientific applications, however, simplified modeling assumptions and missing physics can induce systematic model form errors that cannot be eliminated by uncertainty quantification alone. This mismatch can bias the inferred solution and lead to uncertainty estimates that are overconfident or miscalibrated.
	To address this limitation, we will present an explicit model correction mechanism that
	accounts for the discrepancy between the assumed PDE and the true dynamics.
	
	\section{Latent representation learning based model correction method}\label{Methodology}
	
	In this section, we present LVM-GP with model correction, a probabilistic framework for PDE learning under misspecified physics. The method learns the solution together with a discrepancy term and propagates uncertainty through a shared latent representation. We first describe the model structure, including the shared encoder, the neural operator decoder, and the correction decoder. We then present the learning objective and the inference procedure used to obtain predictions and uncertainty estimates.
	\subsection{Model structure}\label{subsec:model-structure}
	
	We consider the PDE setting in Eq.~\eqref{eq:general_PDE} and assume that the available physics is only approximate, represented by a misspecified differential operator $\widetilde{\mathcal{N}}_{\bm{x}}$. Our goal is to retain this given physics while accounting for the remaining model form discrepancy in a probabilistic and structured manner. The key design is to introduce a \emph{shared} stochastic latent representation: the same latent variable drives the solution model, auxiliary quantities (e.g., boundary terms), and the discrepancy term, so uncertainty is propagated consistently across all quantities of interest.
	
	To this end, we first construct the latent variable through a confidence-aware encoder, following LVM-GP:
	\begin{equation}\label{eq:encoder_z_method}
		\bm{z}_{1}(\bm{x},\bm{\omega}_E; \theta_{E})
		=
		\mathrm{diag}(\bm{m}(\bm{x}; \theta_{m})) \,\bar{\bm{z}}(\bm{x}; \theta_{\bar{\bm{z}}})
		+
		\mathrm{diag}(1-\bm{m}(\bm{x}; \theta_{m})) \,\bm{z}_0(\bm{x},\bm{\omega}_E),
		\qquad
		\theta_{E}=\{\theta_{m},\theta_{\bar{\bm{z}}}\},\quad
		\bm{\omega}_E\sim \mathcal{N}(0,\bm{I}_{M}),
	\end{equation}
	where $\bm{m}(\bm{x})\in[0,1]^M$ is the confidence function, $\bar{\bm{z}}(\bm{x})\in\mathbb{R}^M$ is a deterministic feature map produced by the encoder network, and $\bm{z}_0(\bm{x},\bm{\omega}_E)$ denotes a finite-dimensional approximation to a sample path drawn from a Gaussian process prior with kernel $K(\bm{x},\bm{x}')$. Eq.~\eqref{eq:encoder_z_method} uses $\bm{m}(\bm{x})$ as an input-dependent gate to interpolate between the deterministic feature and the GP component, thereby controlling epistemic uncertainty while preserving GP-induced spatial structure.
	
	Having obtained $\bm{z}_{1}(\bm{x},\bm{\omega}_E)$, we then map it to the solution field using a neural operator decoder. Conditioned on the same latent representation, the decoder produces a conditional Gaussian model for $u$:
	\begin{equation}\label{eq:pred_u}
		\begin{aligned}
			\bm{z}_{1}(\bm{x}, \bm{\omega}_E) &= \bm{z}_{1}(\bm{x},\bm{\omega}_E; \theta_{E}), \\ 
			\bm{z}_{i}(\bm{x}, \bm{\omega}_E)
			&= \sigma\left(
			W_{i}\bm{z}_{i-1}(\bm{x}, \bm{\omega}_E) + b_{i}
			+ \int_{\Omega} k_{i}(\bm{x}, \bm{x'})\bm{z}_{i-1}(\bm{x'},\bm{\omega}_E)\,\mathrm{d}\bm{x'}
			\right), \quad i = 2, \cdots, L-1, \\
			\bm{z}_{L}(\bm{x}, \bm{\omega}_E) &= W_{L}\bm{z}_{L-1}(\bm{x}, \bm{\omega}_E) + b_{L}, \\
			u(\bm{x}, \bm{\omega}_E; \theta_{u})
			&= z_{L}(\bm{x}, \bm{\omega}_E) + \sigma_{u}(\bm{x}; \theta^{u}_{\sigma}) \cdot \varepsilon_{u},
			\qquad
			\varepsilon_{u} \sim \mathcal{N}(0, 1),
		\end{aligned}
	\end{equation}
	where
	\[
	k_{i}(\bm{x}, \bm{x'}) = \frac{\exp(-\alpha_{i}\|\bm{x} - \bm{x'}\|^{2})}{\int_{\Omega}\exp(-\alpha_{i}\|\bm{x} - \bm{y}\|^{2})\,\mathrm{d}\bm{y}}V_{i},
	\qquad
	\theta^{u}_{D} := \{ \alpha_{i}, V_{i}, W_{i}, b_{i} \}^{L-1}_{i = 2} \cup \{ W_{L}, b_{L} \},
	\]
	here $\alpha_i>0$ is a length scale parameter and $V_i$ is a learnable projection matrix, $\theta_{u}=\{\theta_E,\theta_D^u,\theta_\sigma^u\}$. This induces the conditional distribution
	\begin{equation}\label{eq:lvm-gp condition distribution u}
		q^{u}_{D}(u \mid \bm{x}, \bm{\omega}_E; \theta_{u})
		\sim
		\mathcal{N}\Big(\mu_u(\bm{x}, \bm{\omega}_E; \theta_{E}, \theta^{u}_{D}),\ \sigma^2_{u}(\bm{x}; \theta^{u}_{\sigma})\Big),
	\end{equation}
	where $\mu_u(\bm{x},\bm{\omega}_E;\theta_E,\theta_D^u)$ denotes the decoder mean and $\sigma_u^2(\bm{x};\theta_\sigma^u)$ models aleatoric uncertainty.
	
	Along the same latent representation, the boundary term, denoted by $b(\bm{x})$, is modeled in an analogous conditional Gaussian form
	\begin{equation}\label{eq:lvm-gp condition distribution b}
		q^b_{D}(b \mid \bm{x}, \bm{\omega}_E; \theta_b)
		\sim
		\mathcal{N}\Big(\mu_b(\bm{x}, \bm{\omega}_E; \theta_{E}, \theta^{u}_{D}),\ \sigma^2_{b}(\bm{x}; \theta^b_{\sigma})\Big),
		\qquad
		\theta_{b} = \{ \theta_{E}, \theta^{u}_{D}, \theta^{b}_{\sigma} \},
	\end{equation}
	where the mean is induced from the solution mean through the prescribed boundary operator $B_x$, namely
	\begin{equation}\label{eq:mu_b_Bx_mu_u}
		\mu_b(\bm{x}, \bm{\omega}_E; \theta_E, \theta^{u}_{D})
		\coloneqq \mathcal{B}_{\bm{x}}\left[\mu_u(\cdot, \bm{\omega}_E; \theta_E, \theta^{u}_{D})\right](\bm{x}).
	\end{equation}
	
	With the probabilistic model for $(u,b)$ in place, we next incorporate the forcing term through the approximate physics while explicitly accounting for model form discrepancy. Rather than learning the full governing equation end-to-end, we correct the misspecified operator by augmenting $\widetilde{\mathcal{N}}_{x}$ with a learnable correction operator $\mathcal{S}_{\psi}$ that is driven by the same latent representation used to generate the solution:
	\begin{equation}\label{eq:corrected force term}
		\mu_f(\bm{x},\bm{\omega}_E)
		=
		\mathcal{N}_{\bm{x}}[\mu_u(\bm{x}, \bm{\omega}_E)]
		=
		\widetilde{\mathcal{N}}_{\bm{x}}[\mu_u(\bm{x}, \bm{\omega}_E)]
		+
		\mathcal{S}_{\psi}\big(\bm{z}_{1}(\bm{x}, \bm{\omega}_E)\big),
	\end{equation}
	where $\psi$ denotes the parameters of the correction mapping and $\bm{\omega}_E$ is the shared latent randomness introduced by the encoder in Eq.~\eqref{eq:encoder_z_method}. By contrast, $\varepsilon_u$ in Eq.\eqref{eq:pred_u} denotes the independent decoder noise used in the conditional model for $u$.
	
	This construction has two key implications. First, the correction term $\mathcal{S}_{\psi}$ does not introduce an unrelated source of randomness; instead, it is conditioned on the same latent variable $z_1$ as the solution, ensuring that the discrepancy inherits the encoder induced latent structure. Second, because the correction and the solution share $\bm{\omega}_E$, uncertainty in the latent space is propagated coherently to both components, yielding a natural coupling between discrepancy uncertainty and solution uncertainty. As a result, we obtain a corrected conditional Gaussian model for the forcing term:
	\begin{equation}\label{eq:correct_f}
		q^{f}_{D}(f \mid \bm{x}, \bm{\omega}_E; \theta_{f}, \psi)
		\sim
		\mathcal{N}\Big(
		\widetilde{\mathcal{N}}_{\bm{x}}[\mu_u(\bm{x}, \bm{\omega}_E; \theta_{E}, \theta^{u}_{D})] + \mathcal{S}_{\psi}(\bm{z}_{1}(\bm{x}, \bm{\omega}_E)),\ 
		\sigma^{2}_{f}(\bm{x}; \theta^{f}_{\sigma})
		\Big),
		\qquad
		\theta_{f} = \{ \theta_{E}, \theta^{u}_{D}, \theta^{f}_{\sigma} \}.
	\end{equation}
	Notably, the correction term modifies the mean of the operator output, while leaving the variance model $\sigma_f^2(\bm{x};\theta_\sigma^f)$ unchanged, so the uncertainty quantification mechanism remains compatible with the original likelihood design.
	
	\begin{remark}
		The operator level correction in Eq.~\eqref{eq:operator_correction} can be viewed as a special case of our latent-space correction.
		Indeed, if the latent conditioned mapping is restricted to depend on $\bm{z}_1$ only through the decoded solution mean, i.e.,
		\[
		\mathcal{S}_{\psi}\left(\bm{z}_1(\bm{x},\bm{\omega}_E)\right)\equiv 
		\mathcal{C}_{\bm{x}}\left[\mu_u(\cdot,\bm{\omega}_E)\right](\bm{x}),
		\]
		then the corrected forcing mean
		$\mu_f=\widetilde{\mathcal{N}}_{\bm{x}}[\mu_u]+\mathcal{S}_{\psi}(\bm{z}_1)$ reduces to
		$\mu_f=(\widetilde{\mathcal{N}}_{\bm{x}}+\mathcal{C}_{\bm{x}})[\mu_u]$.
	\end{remark}
	
	\begin{remark}
		Note that the correction term is modeled as a pointwise function of the shared latent representation, i.e.,
		$\mathcal{S}_\psi(\bm{x},\bm{\omega}_E)=\mathcal{S}_\psi(\bm{z}_1(\bm{x},\bm{\omega}_E))$.
		This pointwise mapping can be replaced by an integral neural operator structure, analogous to the solution decoder, so that
		$\mathcal{S}_\psi(\bm{z}_1)(\bm{x},\bm{\omega}_E)$ depends on the entire latent field
		$\{\bm{z}_1(\bm{x}',\bm{\omega}_E)\}_{\bm{x}'\in\Omega}$ through an integral kernel:
		\begin{equation}\label{eq:IO_s}
			\begin{aligned}
				\bm{z}^{s}_{i}(\bm{x}, \bm{\omega}_E) &= \sigma \left( W^{s}_{i}\bm{z}^{s}_{i-1}(\bm{x}, \bm{\omega}_{E}) + b^{s}_{i} + \int_{\Omega} k_{i}^{s}(\bm{x}, \bm{x'})\bm{z}^{s}_{i-1}(\bm{x'}, \bm{\omega}_E)\mathrm{d}\bm{x'} \right), \quad i = 2, \cdots, L_2-1, \\
				\bm{z}^{s}_{L_2}(\bm{x}, \bm{\omega}_E) &= W^{s}_{L_2}\bm{z}^{s}_{L_2-1}(\bm{x}, \bm{\omega}_{E}) + b_{L_2}, \\
				\mathcal{S}_\psi(\bm{x}, \bm{\omega}_E) &= \bm{z}^{s}_{L_2}(\bm{x}, \bm{\omega}_E),
			\end{aligned}
		\end{equation}
		with $\bm{z}_1^s(\bm{x},\bm{\omega}_E) = \bm{z}_1(\bm{x},\bm{\omega}_E)$.
		In this case, $\mu_u(\bm{x},\bm{\omega}_E)$ and $\mathcal{S}_\psi(\bm{x},\bm{\omega}_E)$ are generated by two neural operator decoders that share the same latent input field $\bm{z}_1(\bm{x},\bm{\omega}_E)$, but employ different weights and kernels to capture distinct nonlocal couplings for the solution and correction. This separation allows the two outputs to learn different spatial interaction patterns while remaining coupled through the shared latent representation.
		We further include an additional numerical experiment in which the correction decoder adopts the neural operator form in Eq.~\eqref{eq:IO_s}, see Figure~\ref{fig:IO_uands}.
	\end{remark}

	\subsection{Learning and inference}\label{subsec:learning-inference}
	
	For the training process, since all quantities are modeled conditionally on the shared latent randomness, we evaluate the data likelihood via Monte Carlo sampling of $\bm{\omega}_E$, which leads to the empirical log-likelihood terms.
	We define the data objective as
	\[
	\mathcal{L}_{\mathrm{data}}
	=
	w_{u}\mathcal{L}_{\mathrm{data},u}
	+
	w_{f}\mathcal{L}_{\mathrm{data},f}
	+
	w_{b}\mathcal{L}_{\mathrm{data},b},
	\qquad
	\omega_E^{(j)}\sim \mathcal{N}(0,\bm{I}_M),\ j=1,\dots,N_{\bm{\omega}},
	\]
	where $\omega_u, \omega_f,\omega_b$ are prescribed weights, 
	$N_{\bm{\omega}}$ is the number of latent samples, and 
	\begin{align*}
		\mathcal{L}_{\mathrm{data},u}
		&=
		\frac{1}{N_u N_{\bm{\omega}}}\sum_{i=1}^{N_u}\sum_{j=1}^{N_{\bm{\omega}}}
		\log q_D^u\left(\bar u^{(i)}\,\middle|\,\bm{x}_u^{(i)},\bm{\omega}_E^{(j)};\theta_u\right), \\
		\mathcal{L}_{\mathrm{data},b}
		&=
		\frac{1}{N_b N_{\bm{\omega}}}\sum_{i=1}^{N_b}\sum_{j=1}^{N_{\bm{\omega}}}
		\log q_D^b\left(\bar b^{(i)}\,\middle|\,\bm{x}_b^{(i)},\bm{\omega}_E^{(j)};\theta_b\right), \\
		\mathcal{L}_{\mathrm{data},f}
		&=
		\frac{1}{N_f N_{\bm{\omega}}}\sum_{i=1}^{N_f}\sum_{j=1}^{N_{\bm{\omega}}}
		\log q_D^f\left({\bar f}^{(i)}\,\middle|\,\bm{x}_f^{(i)},\bm{\omega}_E^{(j)};\theta_f,\psi\right).
	\end{align*}
	Notably, the forcing likelihood $\mathcal{L}_{\mathrm{data},f}$ merges physical and corrective information by evaluating $\bar f$ against $\widetilde{\mathcal{N}}_{\bm{x}}[\mu_u]+\mathcal{S}_\psi(\bm{z}_1)$, with both parts driven by the shared latent representation.

\begin{figure}[H]
		\centering
		\includegraphics[width=1.0\textwidth]{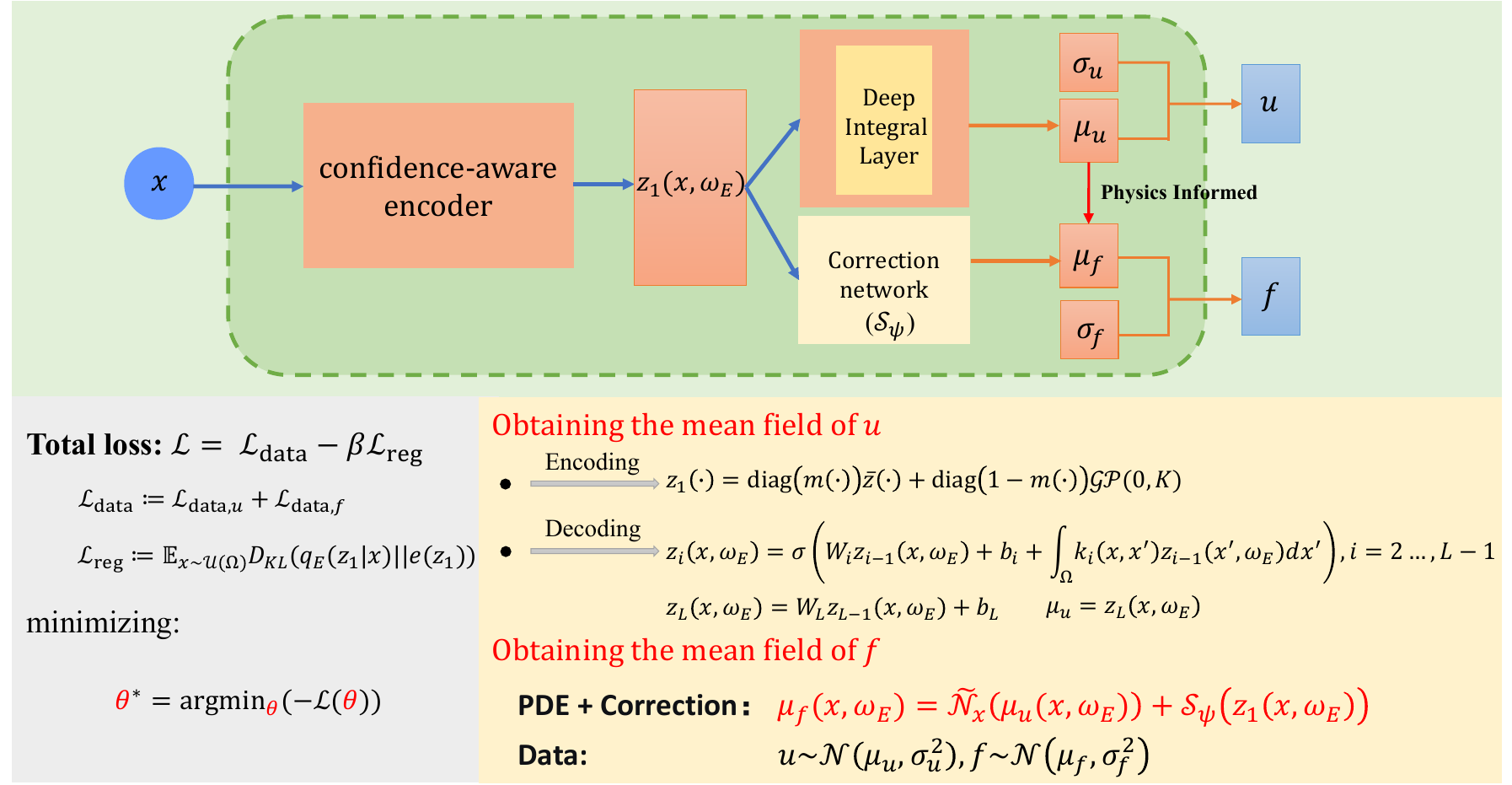} 
		\caption{Schematic of latent-space model correction in the LVM-GP framework. The input $\bm{x}$ is encoded by a confidence-aware encoder to form a stochastic latent code $\bm{z}_{1}(\bm{x},\bm{\omega})$, which is shared by two decoders.
			The solution decoder is implemented by stacked deep integral layers, producing the predictive mean and uncertainty $(\mu_u,\sigma_u)$ of the solution field $u(\bm{x})$.
			In parallel, the correction decoder outputs a discrepancy term via the correction network $\mathcal{S}_{\psi}$, and the physics quantity is coupled to the solution through the corrected operator relation
			$\mu_f=\widetilde{\mathcal{N}}_{\bm{x}}(\mu_u)+\mathcal{S}_{\psi}(\bm{z}_{1}(\bm{x},\bm{\omega}))$.
			Finally, by leveraging learned variance functions, the model outputs complete probability distributions for the solution $u(\bm{x})$, the source term $f(\bm{x})$, and the correction term $\mathcal{S}_{\psi}(\bm{z}_{1}(\bm{x}, \bm{\omega}))$, quantifying their uncertainties.}
		\label{fig: clvm-gp}
	\end{figure}

	To better regularize the latent representation and avoid degeneracy and overfitting, we introduce a latent regularization term on the encoder distribution. More specifically,
	note that the encoder in Eq.~\eqref{eq:encoder_z_method} follows
	\[
	q_E(\bm{z}_1|\bm{x})=\mathcal{N}\left(\mathrm{diag}(\bm{m}(\bm{x}))\,\bar {\bm{z}}(\bm{x}),\ \mathrm{diag}(1-\bm{m}(\bm{x}))^2\right).
	\]
	Thus we can compute the KL divergence between $q_E(\bm{z}_1|\bm{x})$ and marginal distribution $e(\bm{z}_1)$ on the entire
	physical domain as a regularization term, namely,\[
	\mathcal{L}_{\mathrm{reg}}
	=
	\mathbb{E}_{\bm{x}\sim\mathcal{U}(\Omega)}\,D_{\mathrm{KL}}\left[q_E(\bm{z}_1|\bm{x})\,\|\,e(\bm{z}_1)\right],
	\qquad
	e(\bm{z}_1)=\mathcal{N}(0,\bm{I}_M),
	\]
	and approximate it empirically by $\bm{x}_i\sim\mathcal{U}(\Omega)$ and $\bm{z}_{1,i}\sim q_E(\cdot|\bm{x}_i)$:
	\[
	\mathcal{L}_{\mathrm{reg}}
	\approx
	\frac{1}{N}\sum_{i=1}^{N}\Big[\log q_E(z_{1,i}|x_i)-\log e_1(z_{1,i})\Big].
	\]
	Finally, we maximize the overall objective
	\[
	\mathcal{L}=\mathcal{L}_{\mathrm{data}}-\beta\,\mathcal{L}_{\mathrm{reg}}.
	\]
	An overview of the proposed latent representation learning based model correction framework is illustrated in Figure~\ref{fig: clvm-gp}, and the corresponding training procedure is summarized in Algorithm~\ref{algorithm}.
	\begin{algorithm}[H]
		\begin{itemize}
			\item\textbf{1. Input:} Training set $\{\bm{x}^{(i)}_{u}, \bar{u}^{(i)} \}^{N_{u}}_{i = 1}$, $ \{ \bm{x}^{(i)}_{f}, \bar{f}^{(i)} \}^{N_{f}}_{i = 1}$, $\{ \bm{x}^{(i)}_{b}, \bar{b}^{(i)} \}^{N_{b}}_{i = 1}$, learning rate $\eta$, and regularization parameter $\beta$
			\item\textbf{2. While not converged do}
			\item\textbf{3.}\quad Draw samples from $\bm{\omega^}{(j)}_{E} \sim \mathcal{N}(0, \bm{I_{M}}), j = 1, \dots, N_{\bm{\omega}}$
			\item\textbf{4.}\quad Compute the total loss: 
			$$
			\mathcal{L} = \omega_{u}\cdot\mathcal{L}_{\mathrm{data},u} + \omega_{f}\cdot\mathcal{L}_{\mathrm{data},f} + \omega_{b}\cdot\mathcal{L}_{\mathrm{data},b} - \beta \cdot \mathcal{L}_{\mathrm{reg}}
			$$
			where $\mathcal{L}_{\mathrm{data},f}$ is computed using the corrected distribution (see Eq.~\eqref{eq:correct_f})
			\item\textbf{5.}\quad Update parameters $\theta_{u}, \theta_{f}, \theta_{b}, \psi$ via gradient descent:
			$$
			W \leftarrow W + \eta \cdot \frac{\partial \mathcal{L}}{\partial W}
			$$
			\item\textbf{6. End While}
		\end{itemize}
		\caption{Latent representation learning based model correction}
		\label{algorithm}
	\end{algorithm}
	
	After training, we perform inference via Monte Carlo propagation through the shared latent space.
	We draw $\bm{\omega}_E^{(j)}\sim\mathcal{N}(0,\bm{I}_M)$ for $j=1,\dots,N_{\bm{\omega}}$ and generate
	\[
	s^{(j)}(\bm{x}) \coloneqq \mathcal{S}_{\psi}\left(\bm{z}_{1}(\bm{x},\bm{\omega}_E^{(j)})\right),\qquad
	u^{(j)}(\bm{x}) \coloneqq \mu_u(\bm{x},\bm{\omega}_E^{(j)})+\sigma_u(\bm{x})\,\varepsilon_u^{(j)},\quad
	\varepsilon_u^{(j)}\sim\mathcal{N}(0,1),
	\]
	and similarly for other quantities. For any quantity $g\in\{u,\mathcal{S}_\psi\}$, we estimate its predictive mean and variance by
	\[
	\widehat{\mathbb{E}}[g](\bm{x})=\frac{1}{N_{\bm{\omega}}}\sum_{j=1}^{N_{\bm{\omega}}} g^{(j)}(\bm{x}),\qquad
	\widehat{\mathrm{Var}}[g](\bm{x})=\frac{1}{N_{\bm{\omega}}}\sum_{j=1}^{N_{\bm{\omega}}}
	\Big(g^{(j)}(\bm{x})-\widehat{\mathbb{E}}[g](\bm{x})\Big)^2.
	\]
	Here $\sigma_u(\bm{x})$ accounts for aleatoric uncertainty in $u$, while $\mathcal{S}_{\psi}$ reflects only latent-driven epistemic uncertainty.

	\begin{remark}
		When observation noise is negligible, we set the aleatoric models to zero or remove them, i.e.,
		$\sigma_u(\bm{x})\equiv 0$ and $\sigma_f(\bm{x})\equiv 0$, and similarly for $b$ if applicable.
		In this case, the likelihood terms reduce to MSE-type objectives, and the predictive sampling simplifies to
		$u^{(j)}(\bm{x})=\mu_u(\bm{x},\bm{\omega}_E^{(j)})$ and
		$f^{(j)}(\bm{x})=\widetilde{\mathcal{N}}_{\bm{x}}[\mu_u(\cdot,\bm{\omega}_E^{(j)})](\bm{x})+s^{(j)}(\bm{x})$.
		Uncertainty is then purely epistemic and is quantified by the dispersion across latent samples
		$\{\bm{\omega}_E^{(j)}\}_{j=1}^{N_{\bm{\omega}}}$.
	\end{remark}
	
	\section{Numerical experiments}\label{Experiments and Results}
	In this section, we evaluate the proposed model correction method based on latent representation learning within the LVM-GP framework on four representative problems, including ordinary differential equations, reaction-diffusion equations, two-dimensional channel flows, and two-dimensional cavity flows. The experiments are designed to assess prediction accuracy, uncertainty quantification, and model correction performance under noisy data and model misspecification.
	\subsection{An ODE System}\label{ODE_System}
	We first validate the proposed method on a one-dimensional ordinary differential equation (ODE):
	\begin{equation}\label{eq:example_ode}
		\frac{du}{dt} = f(t) + \lambda u(1-u), \quad t\in[0,1], \qquad u(0)=u_0,
	\end{equation}
	where $\lambda \ge 0$ is a constant parameter (ground truth $\lambda^*=1.5$), $u_0$ is the initial condition, and $f(t)$ is a prescribed source term.
	In practice, simplified modeling assumptions and limited observations can lead to discrepancies between the true dynamics and the assumed model.
	To study the impact of such model form errors, we intentionally introduce misspecification in the reaction term and consider the following three settings:
	\begin{itemize}
		\item S1: The reaction term is correctly specified as $\lambda u(1-u)$, and $\lambda\ge 0$ is treated as unknown.
		\item S2: The reaction term is misspecified as $\lambda \cos(u)$, and $\lambda\ge 0$ is treated as unknown.
		\item S3: The reaction term is misspecified as $\lambda \cos(u)$ with $\lambda$ set to $0.2$. To compensate for this misspecification, we introduce a correction term $\mathcal{S}_{\psi}$, resulting in $0.2\cos(u)+\mathcal{S}_{\psi}$.
	\end{itemize}
	
	To simplify the notations, we write the reaction term in Eq.~\eqref{eq:example_ode} as $\phi$. 
	In S1 and S2, we evaluate $\phi$ using the LVM-GP estimate $\tilde{\lambda}$ (see \cite{feng2025lvm}) and the learned solution $u_\theta$, namely
	$\phi=\tilde{\lambda}\,u_{\theta}(1-u_{\theta})$ and $\phi=\tilde{\lambda}\cos(u_{\theta})$, respectively. In S3, we keep the misspecified coefficient $\lambda=0.2$ and augment the reaction term with a latent-conditioned correction:
	\[
	\phi=0.2\cos(u_{\theta})+\mathcal{S}_{\psi}\left(\bm{z}_{1})\right).
	\]
	Here $\mathcal{S}_{\psi}$ accounts for the discrepancy introduced by the misspecification.	  
	
	The encoder and the solution decoder are implemented as fully connected neural networks, each comprising three hidden layers with 64 neurons and using the Mish activation function. The latent dimension is set to 20. The correction decoder is implemented by another fully connected network with two hidden layers of 64 neurons and the Tanh activation function. Unless otherwise specified, we use these network architectures in all experiments.
	During training, we employ the Adam optimizer with an initial learning rate of $0.001$, which decays by a factor of $0.7$ every 1,000 steps, and we run 20,000 training iterations in all experiments.
	For noise free experiments, we set $\beta=1\times 10^{-5}$ and optimize the predictive mean only throughout the entire training process.
	For noisy data experiments, we set $\beta=0.05$ and adopt a two-phase schedule: the first 10,000 iterations optimize the predictive mean only, and the subsequent 10,000 iterations jointly optimize the predictive mean and the predictive standard deviation.
	
\begin{figure}[H]
		\begin{center}
			\begin{overpic}[width=0.33\textwidth, trim=0 0 0 0, clip=True]{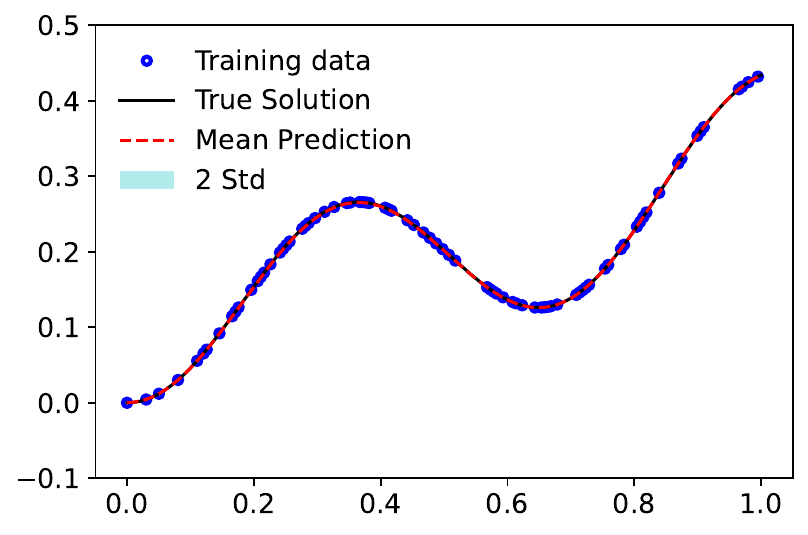}
				\put(52.5,67.5){{$u$}}
			\end{overpic}
			\begin{overpic}[width=0.33\textwidth, trim=0 0 0 0, clip=True]{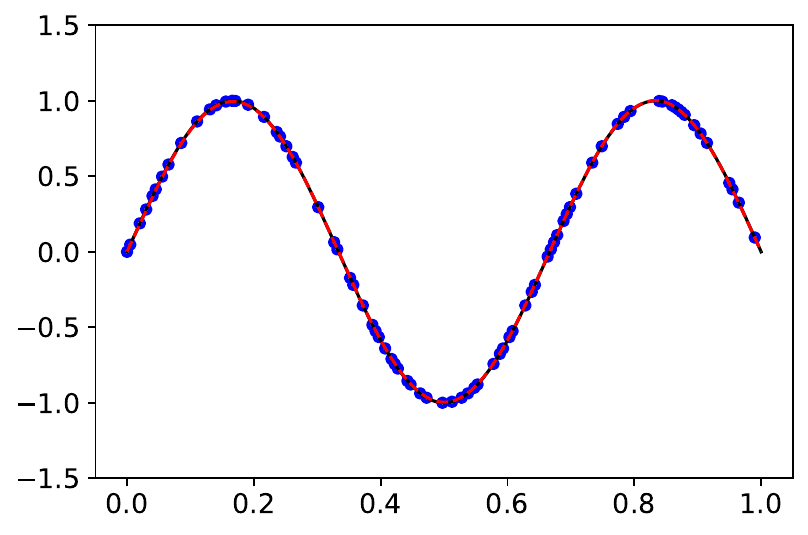}
				\put(52.5,67.5){{$f$}}
			\end{overpic}
			\begin{overpic}[width=0.33\textwidth, trim=0 0 0 0, clip=True]{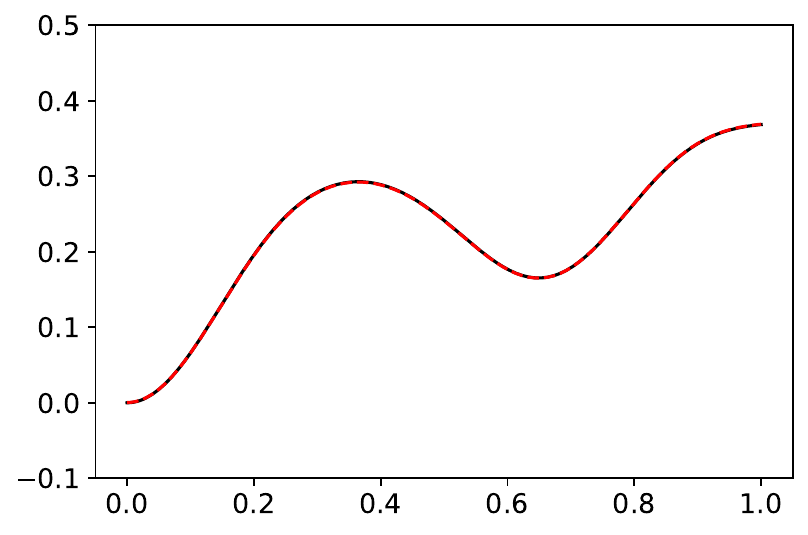}
				\put(50,67.5){{$\phi$}}
			\end{overpic}
			\subcaption{S1: LVM-GP with the correct physical model}
			
			\begin{overpic}[width=0.33\textwidth, trim=0 0 0 0, clip=True]{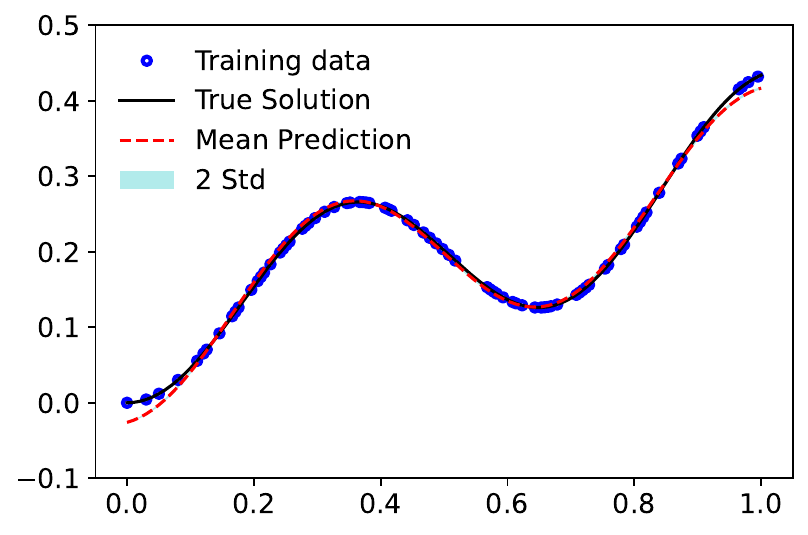}
			\end{overpic}
			\begin{overpic}[width=0.33\textwidth, trim=0 0 0 0, clip=True]{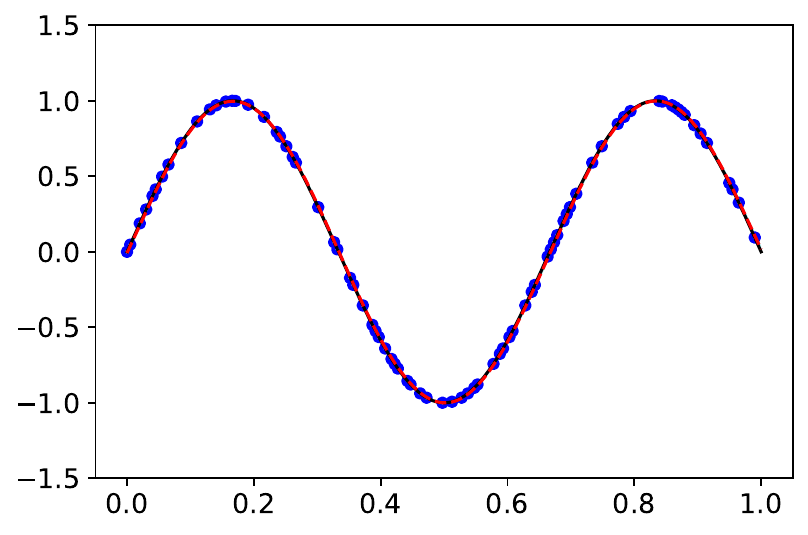}
			\end{overpic}
			\begin{overpic}[width=0.33\textwidth, trim=0 0 0 0, clip=True]{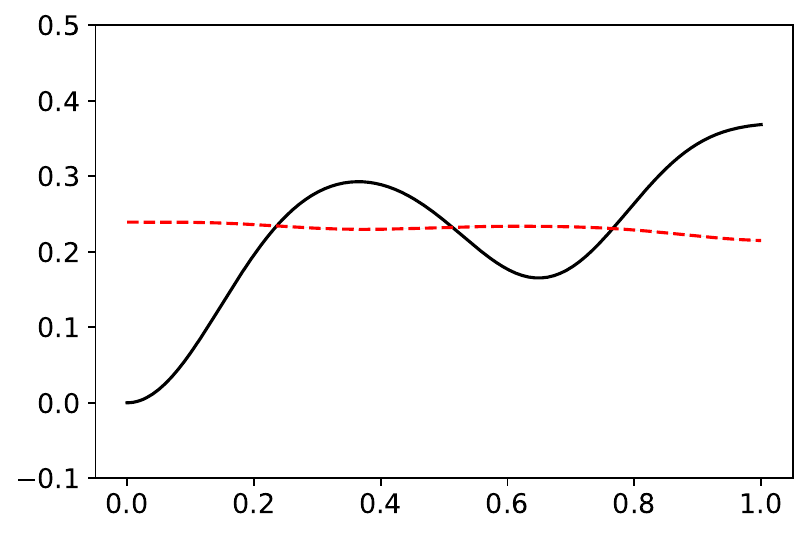}
			\end{overpic}
			\subcaption{S2: LVM-GP with a misspecified physical model}
			
			\begin{overpic}[width=0.33\textwidth, trim=0 0 0 0, clip=True]{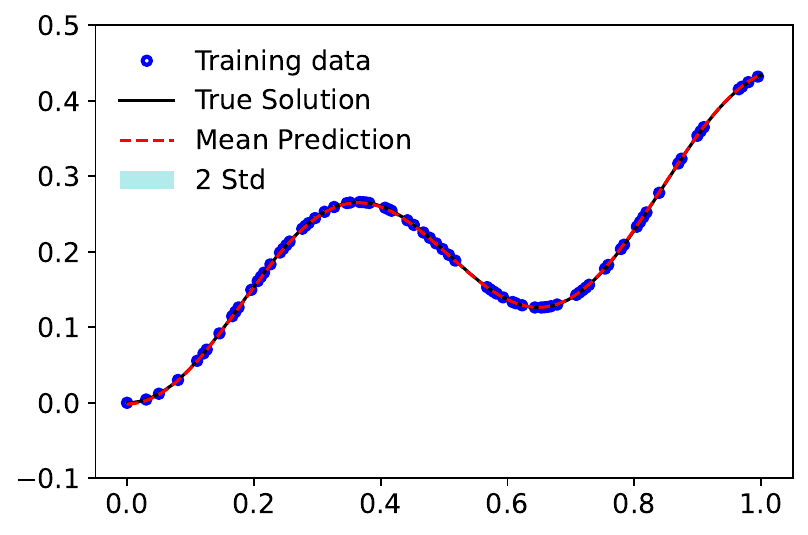}
			\end{overpic}
			\begin{overpic}[width=0.33\textwidth, trim=0 0 0 0, clip=True]{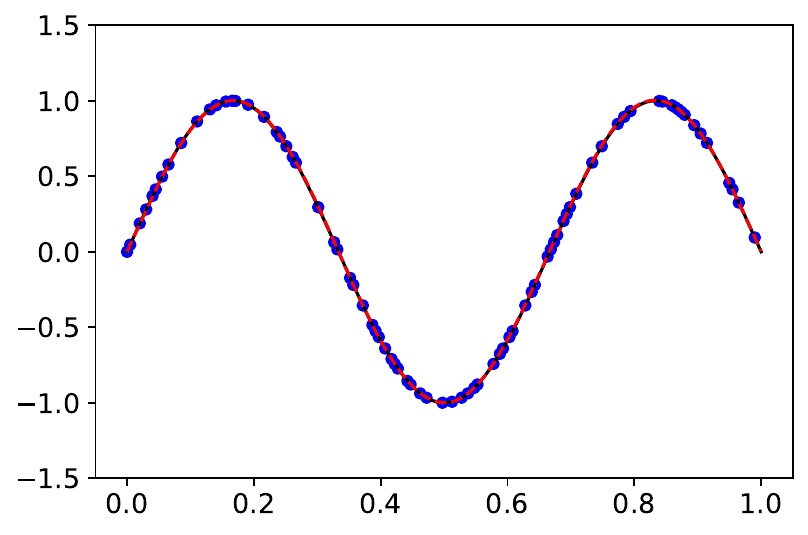}
			\end{overpic}
			\begin{overpic}[width=0.33\textwidth, trim=0 0 0 0, clip=True]{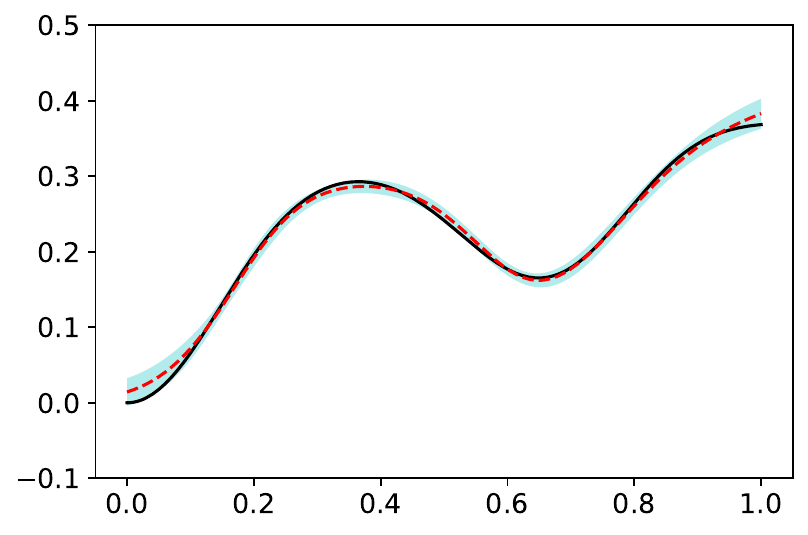}
			\end{overpic}
			\subcaption{S3: Correcting the misspecified physical model by a correction network using LVM-GP}
		\end{center}
		\caption{ODE system. Predictions of $u$, $f$, and $\phi$ under three modeling scenarios using the same noise-free dataset: (a) LVM-GP with the correct physical model, (b) LVM-GP with a misspecified physical model, and (c) correcting the misspecified physical model by a correction network using LVM-GP. Blue circles represent training data, the solid black line denotes the true solution, the red dashed line indicates the model’s mean prediction, and the shaded region corresponds to the ±2 standard deviation predictive uncertainty interval.}
		\label{ode_noisefree}
	\end{figure}

	We begin with a noise free, data sufficient setting to provide a reference benchmark. The dataset is generated by solving Eq.~\eqref{eq:example_ode} with MATLAB's \texttt{ode45}~\cite{inc2022matlab}, where $f(t)=\sin(3\pi t)$ and $u_0=0$, and we take 80 randomly sampled observations of $u(t)$ and $f(t)$ on $[0,1]$. Under the correct reaction model (S1), LVM-GP accurately infers $\phi$ and produces a surrogate solution $u_\theta$ that is consistent with both the observations and the governing ODE (Figure~\ref{ode_noisefree}(a)). When the reaction term is misspecified (S2), 
	the inferred $\phi$ becomes inaccurate, leading to a clear degradation in performance (Figure~\ref{ode_noisefree}(b)). 
	S3 mitigates this mismatch by augmenting the misspecified reaction term with a learned correction network. 
	This leads to a clear improvement in the estimation of $\phi$ (\ref{ode_noisefree}(c)). In addition, Table~\ref{tab:ode_full_noisefree} shows that S3 attains prediction accuracies for $u$ and $f$ that are comparable to those of the correct-model baseline S1. Further analysis for clean but sparse observations is provided in \ref{Additional_results_for_ode_system}.

	\begin{table}[H]
		\centering
		\begin{tabular}{c|c|c|c|c|c}
			\hline
			\hline
			& $\tilde{\lambda}$ & Error of $\phi$ & Error of $u$ & Error of $f$ & Error of $\tilde{u}$ \\
			\hline
			S1: Known model & 1.5004 $\pm$ 0.0007 & 0.0015 & 0.0016 & 0.0040 & 0.0084 \\
			\hline
			S2: Misspecified model & 0.2370 $\pm$ 0.0017 & 0.4182 & 0.0316 & 0.0110 & 0.1072 \\
			\hline
			S3: Corrected model & 0.2 & 0.0274 & 0.0034 & 0.0050 & 0.0066 \\
			\hline
			\hline
		\end{tabular}
		\caption{ODE system. Predictions of $\phi, u, f$ and the reconstructed state $\tilde{u}$ are presented, using a noise-free dataset. The table presents the relative $L_2$ errors between the reference solution and the mean predictions from LVM-GP (S1 and S2) as well as from the LVM-GP based model correction (S3). The reconstructed state $\tilde{u}$ is obtained by solving the identified ODE system using the learned $\phi$, the inferred $f$, and the initial condition $u_0$, with $f$ represented continuously via cubic spline interpolation of its discrete model predictions.}
		\label{tab:ode_full_noisefree}
	\end{table}
	
\begin{figure}[H]
		\begin{center}
			\begin{overpic}[width=0.33\textwidth, trim=0 0 0 0, clip=True]{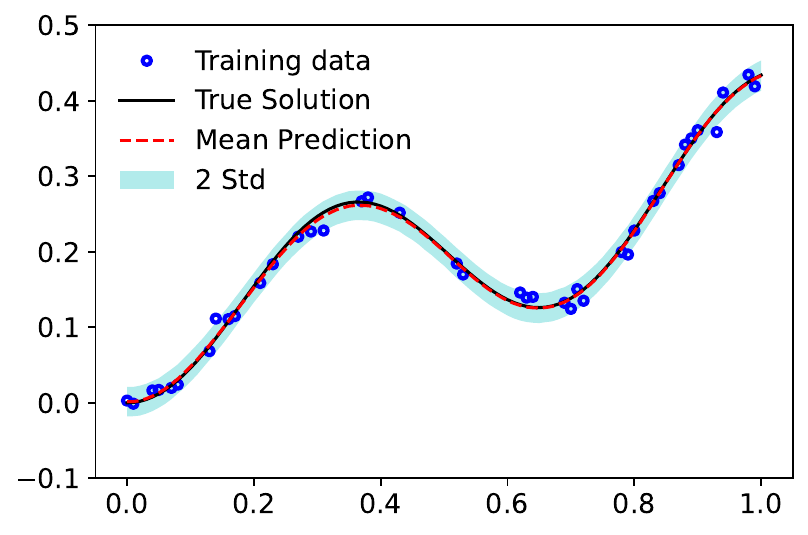}
				\put(52.5,67.5){{$u$}}
			\end{overpic}
			\begin{overpic}[width=0.33\textwidth, trim=0 0 0 0, clip=True]{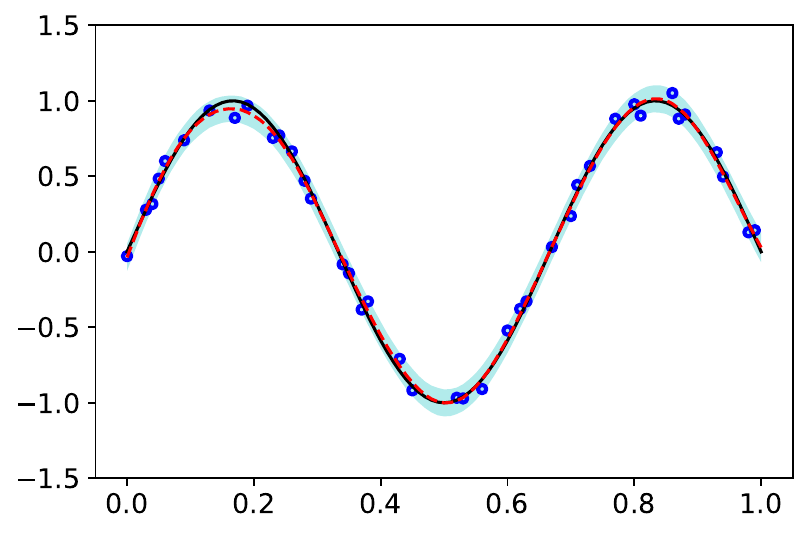}
				\put(52.5,67.5){{$f$}}
			\end{overpic}
			\begin{overpic}[width=0.33\textwidth, trim=0 0 0 0, clip=True]{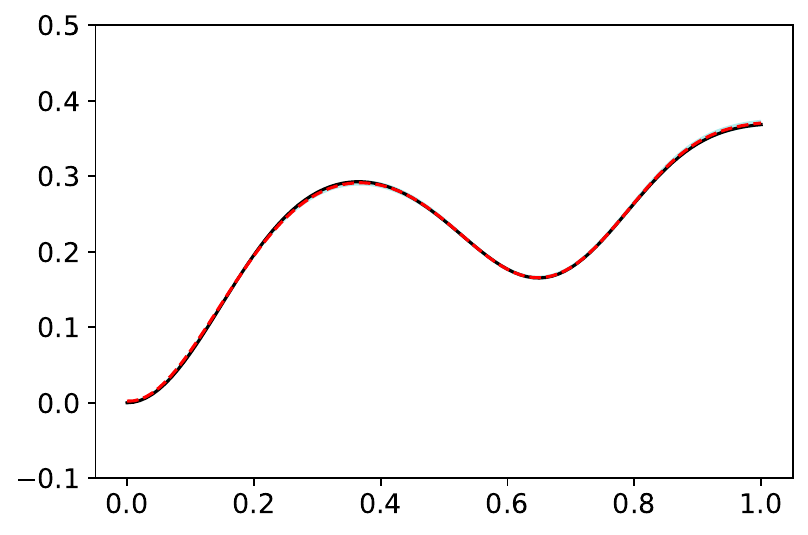}
				\put(50,67.5){{$\phi$}}
			\end{overpic}
			\subcaption{S1: LVM-GP with the correct physical model}
			
			\begin{overpic}[width=0.33\textwidth, trim=0 0 0 0, clip=True]{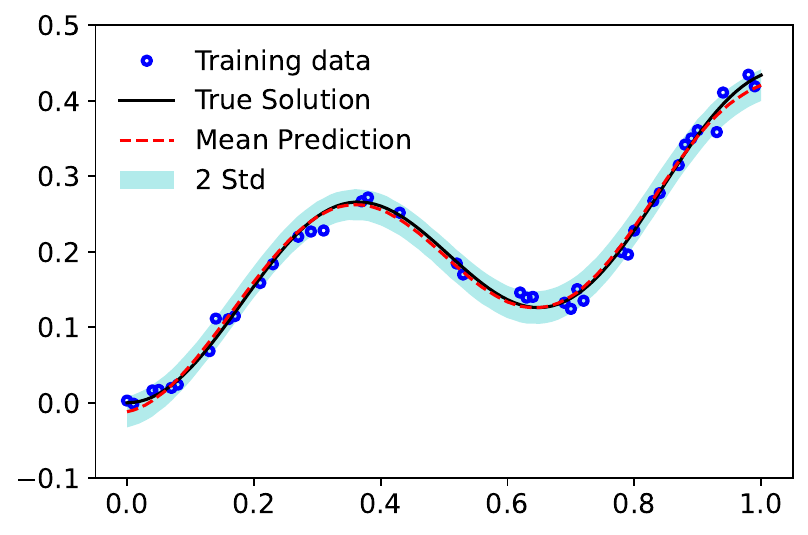}
			\end{overpic}
			\begin{overpic}[width=0.33\textwidth, trim=0 0 0 0, clip=True]{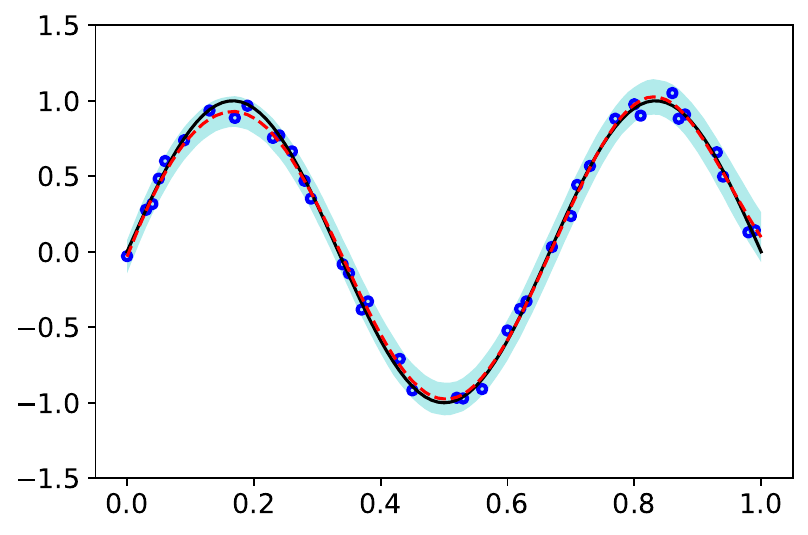}
			\end{overpic}
			\begin{overpic}[width=0.33\textwidth, trim=0 0 0 0, clip=True]{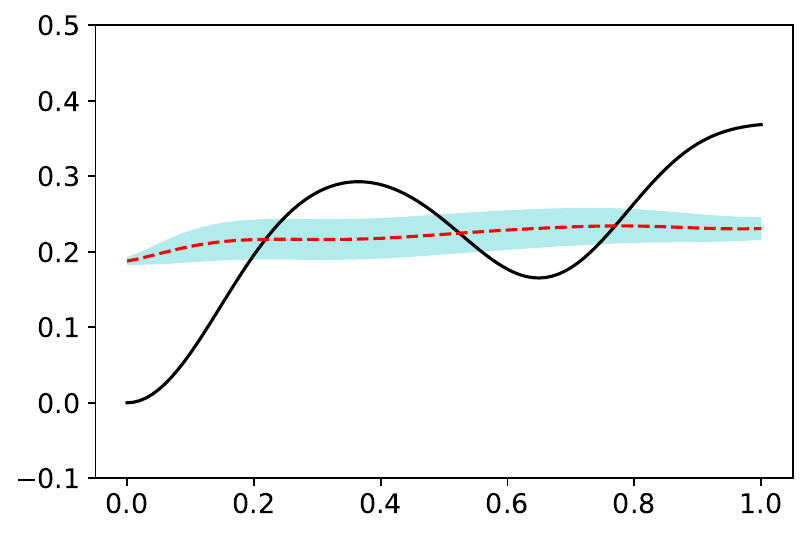}
			\end{overpic}
			\subcaption{S2: LVM-GP with a misspecified physical model}
			
			\begin{overpic}[width=0.33\textwidth, trim=0 0 0 0, clip=True]{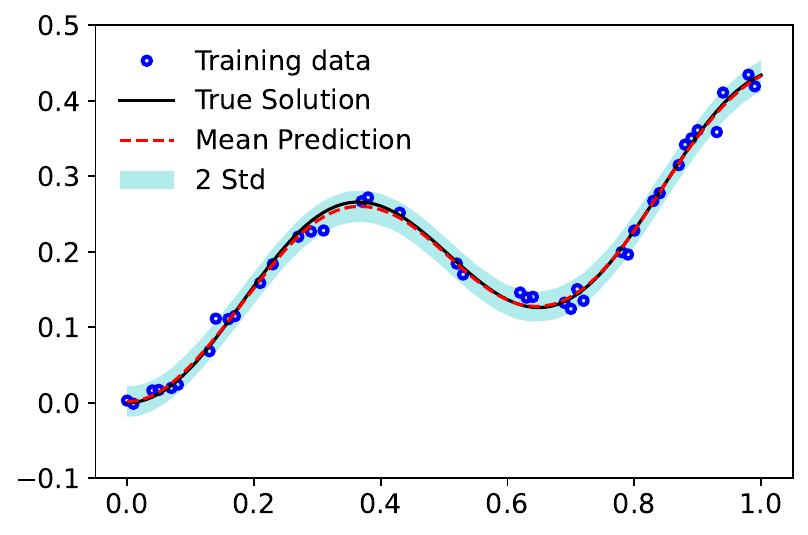}
			\end{overpic}
			\begin{overpic}[width=0.33\textwidth, trim=0 0 0 0, clip=True]{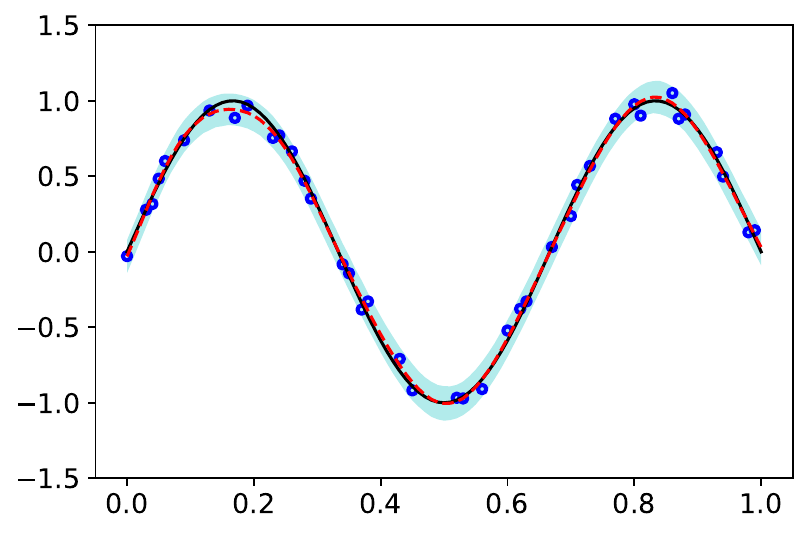}
			\end{overpic}
			\begin{overpic}[width=0.33\textwidth, trim=0 0 0 0, clip=True]{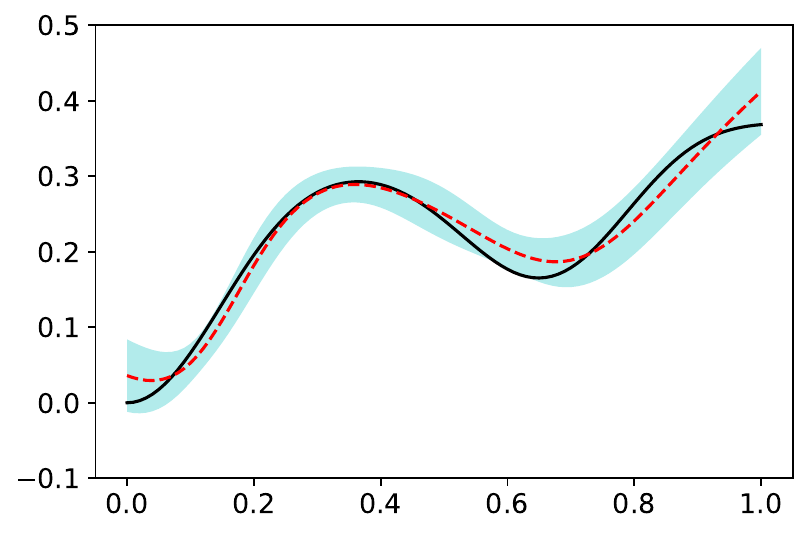}
			\end{overpic}
			\subcaption{S3: Correcting the misspecified physical model by a correction network using LVM-GP}

			\begin{overpic}[width=0.33\textwidth, trim=0 0 0 0, clip=True]{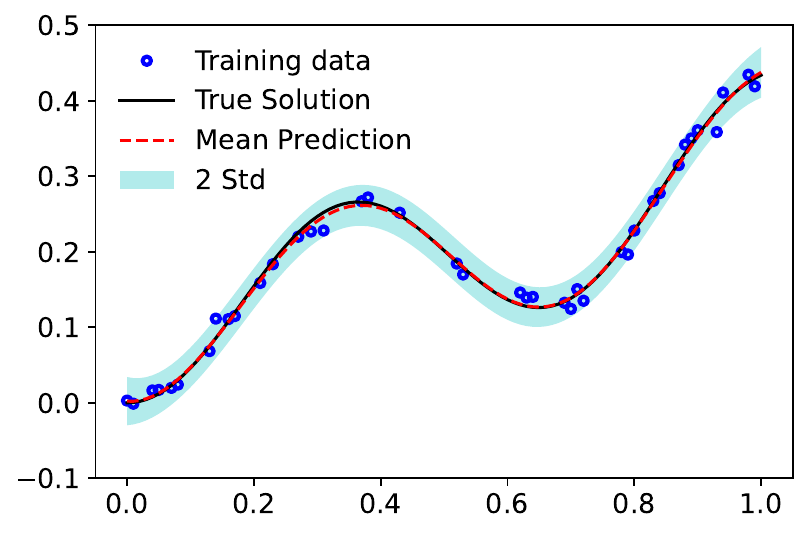}
			\end{overpic}
			\begin{overpic}[width=0.33\textwidth, trim=0 0 0 0, clip=True]{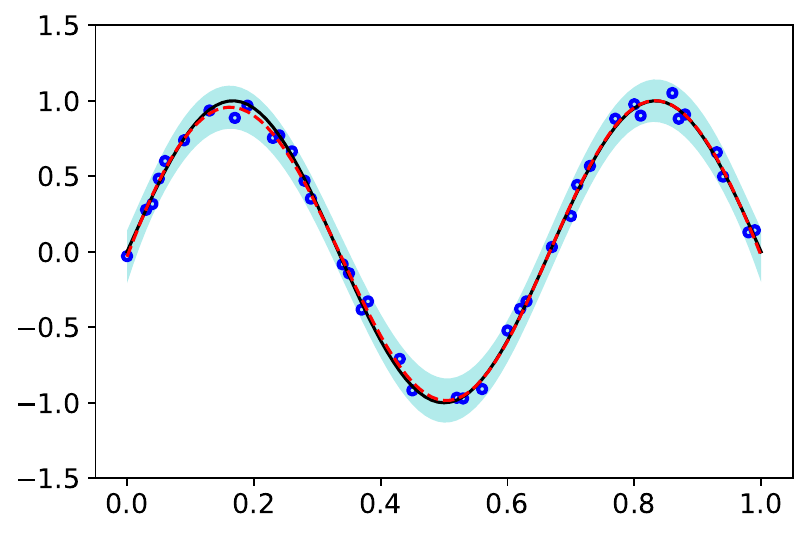}
			\end{overpic}
			\begin{overpic}[width=0.33\textwidth, trim=0 0 0 0, clip=True]{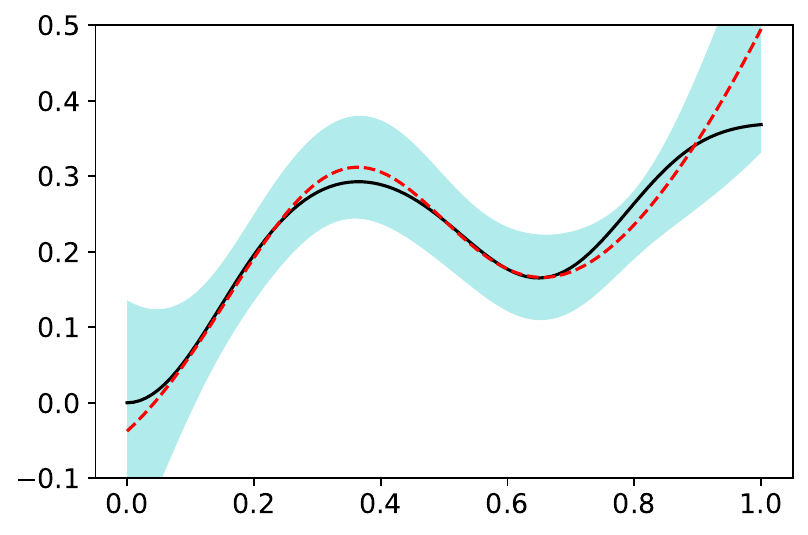}
			\end{overpic}
			\subcaption{S4: Correcting the misspecified physical model by a correction network using B-PINNs}

		\end{center}
		\caption{ODE system. Predictions of $u$, $f$, and $\phi$ under three modeling scenarios using the same sparse and noisy dataset: (a) LVM-GP with the correct physical model, (b) LVM-GP with a misspecified physical model, (c) correcting the misspecified physical model by a correction network using LVM-GP and (d) correcting the misspecified physical model by a correction network using B-PINNs. Blue circles represent training data, the solid black line denotes the true solution, the red dashed line indicates the model’s mean prediction, and the shaded region corresponds to the ±2 standard deviation predictive uncertainty interval.}
		\label{ode_noise}
	\end{figure}

	We next consider sparse and noisy observations while keeping all other settings unchanged. Over the interval $t\in[0,1]$, we randomly sample 40 observations of $u$ and $f$, and perturb each measurement with independent additive Gaussian noise with known standard deviations ($\sigma=0.01$ for $u$ and $\sigma=0.05$ for $f$). The corresponding results are reported in Figure~\ref{ode_noise}. With the correct physical model (S1), the method reconstructs $u$, $f$, and $\phi$ accurately, and the pointwise errors are consistently covered by the 
	predicted uncertainty bands (Figure~\ref{ode_noise}(a)). In contrast, under the misspecified model (S2), the predictions degrade and the associated uncertainty estimates become unreliable (Figure~\ref{ode_noise}(b)). When the misspecified reaction term is augmented using the proposed model-correction framework (S3), the prediction accuracy improves substantially, and the errors again remain within the estimated uncertainty range (Figure~\ref{ode_noise}(c)).Compared with S1, the uncertainty in $\phi$ is consistently larger in S3, which is reasonable since it reflects the additional epistemic uncertainty induced by the lack of knowledge about the reaction model. Figure~\ref{ode_noise}(d) further demonstrates the performance of B-PINNs in correcting a misspecified physical model (S4). B-PINNs not only provide accurate predictions of the solution and effectively quantify its uncertainty, but also yield accurate estimates of the correction term together with reliable robustness measures. Our method achieves comparable performance in correcting model discrepancies, with similar prediction accuracy and uncertainty quantification capability, while avoiding the need for posterior sampling.
	
\begin{figure}[H]
		\centering
		\begin{subfigure}[b]{0.3\textwidth}
			\centering
			\includegraphics[width=\textwidth, trim=0 0 0 0, clip]{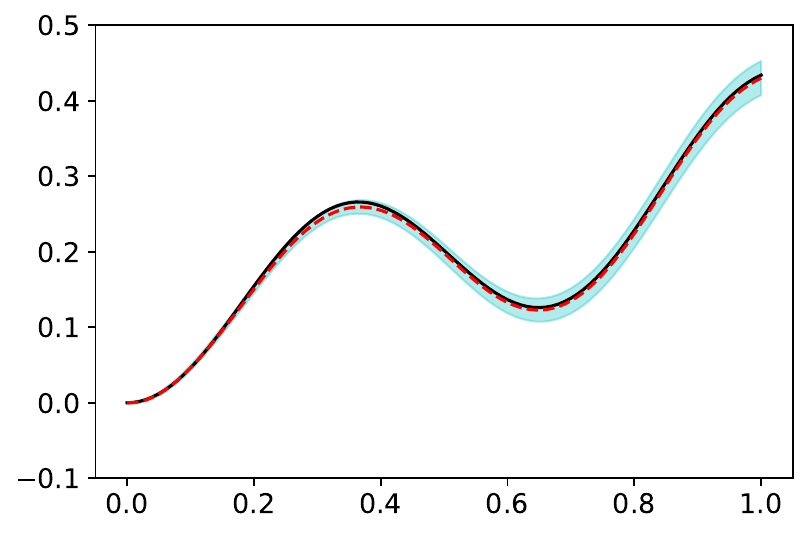}
			\put(-155, 45) {$\tilde{u}$}
			\caption{S1}
		\end{subfigure}
		\hfill 
		\begin{subfigure}[b]{0.3\textwidth}
			\centering
			\includegraphics[width=\textwidth, trim=0 0 0 0, clip]{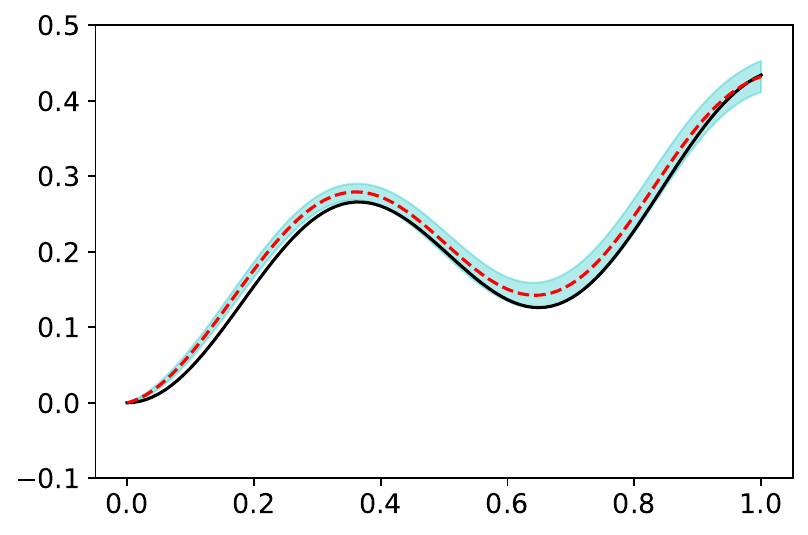}
			\caption{S2}
		\end{subfigure}
		\hfill
		\begin{subfigure}[b]{0.3\textwidth}
			\centering
			\includegraphics[width=\textwidth, trim=0 0 0 0, clip]{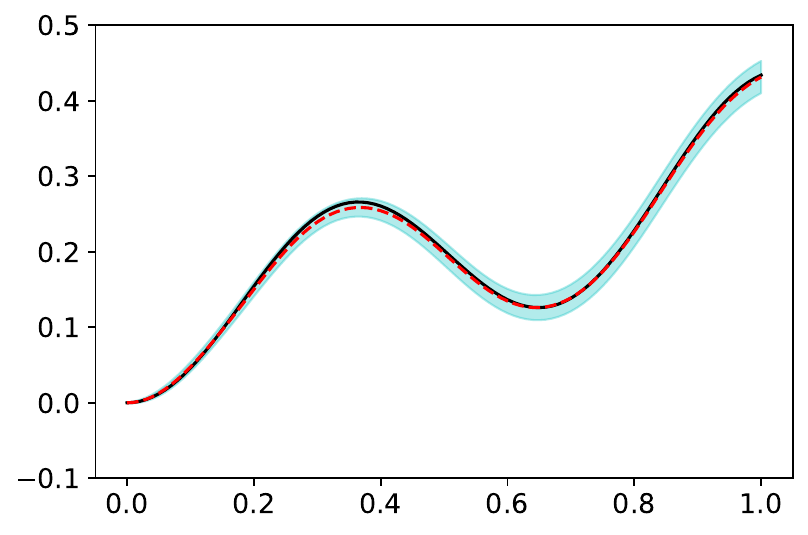}
			\caption{S3}
		\end{subfigure}
		
		\caption{\textbf{ODE system}. Predictions of $\tilde{u}$ under different scenarios using the same sparse and noisy dataset. (a) LVM-GP with the correct physical model. (b) LVM-GP with a misspecified physical model. (c) LVM-GP with a correction network that compensates for model misspecification. The true solution is shown as a solid black line, the model's mean prediction as a red dashed line, and the predictive uncertainty interval (±2 standard deviations) as a shaded region.}
		\label{fig:ode_noise_tildeu}
	\end{figure}

	To validate the correctness of the model inferred from the data, we reconstructed 
	$\tilde{u}$ following the same procedure described in \cite{zou2024correcting}. As shown in Figure~\ref{fig:ode_noise_tildeu}, the reconstruction of $\tilde{u}$ is based on 5,000 samples of $f$ and $\phi$ obtained from the LVM-GP (for S1 and S2) or the model-corrected version based on LVM-GP (for S3). The uncertainty in the reconstructed $\tilde{u}$ is directly propagated from the uncertainty in the predicted $f$ and $\phi$. The relative $L_2$ errors between the mean of $\tilde{u}$ from model prediction and the reference solution are also summarized in Table~\ref{tab:ode_gappy_noise}. 
	As expected, LVM-GP with the correctly specified physics attains the smallest reconstruction error, while the misspecified model performs the worst. Introducing model correction substantially reduces the error and brings the reconstruction accuracy close to that of the correct-model baseline.  	
	
	\begin{table}[H]
		\centering
		\begin{tabular}{c|c|c|c|c|c}
			\hline
			\hline
			& $\tilde{\lambda}$ & Error of $\phi$ & Error of $u$ & Error of $f$ & Error of $\tilde{u}$ \\
			\hline
			S1: Known model & 1.5095 $\pm$ 0.0087 & 0.0063 & 0.0117 & 0.0313 & 0.0158 \\
			\hline
			S2: Misspecified model & 0.2305 $\pm$ 0.0196 & 0.3632 & 0.0236 & 0.0486 & 0.0725 \\
			\hline
			S3: Corrected model using our approach & 0.2 & 0.0723 &  0.0160 & 0.0341
			& 0.0180 \\
			\hline
			S4: Corrected model using B-PINNs & 0.2 & 0.1114 &  0.0112 & 0.0289
			& 0.0173 \\
			\hline
			\hline
		\end{tabular}
		\caption{ODE system. Predictions of $\phi, u, f$ and the reconstructed state $\tilde{u}$ with sparse and noisy data. The table presents the relative $L_2$ errors between the reference solution and the mean predictions from LVM-GP (S1 and S2) as well as from the LVM-GP-based model correction (S3). The reconstructed state $\tilde{u}$ is obtained by solving the identified ODE system using the learned $\phi$, the inferred $f$, and the initial condition $u_0$, with $f$ represented continuously via cubic spline interpolation of its discrete model predictions.}
		\label{tab:ode_gappy_noise}
	\end{table}
	
\begin{figure}[H]
		\centering
		\begin{subfigure}[b]{0.3\textwidth}
			\centering
			\includegraphics[width=\textwidth, trim=0 0 0 0, clip]{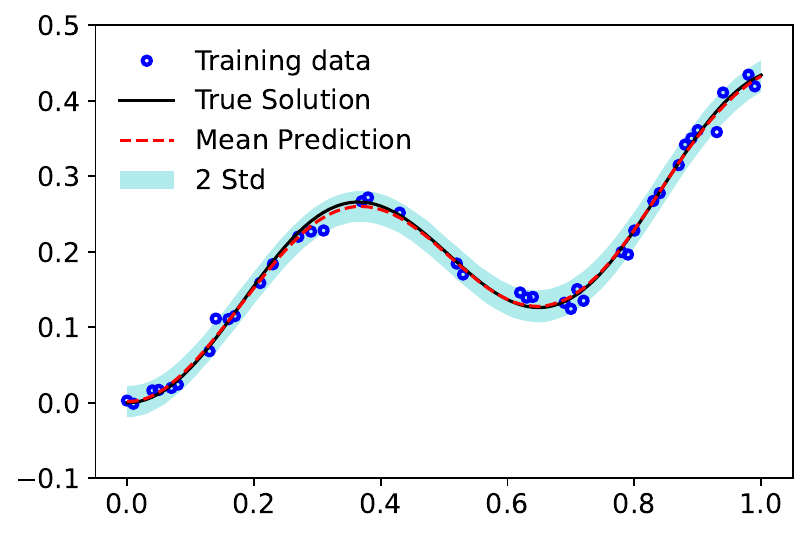}
			\caption{$u$}
		\end{subfigure}
		\hfill 
		\begin{subfigure}[b]{0.3\textwidth}
			\centering
			\includegraphics[width=\textwidth, trim=0 0 0 0, clip]{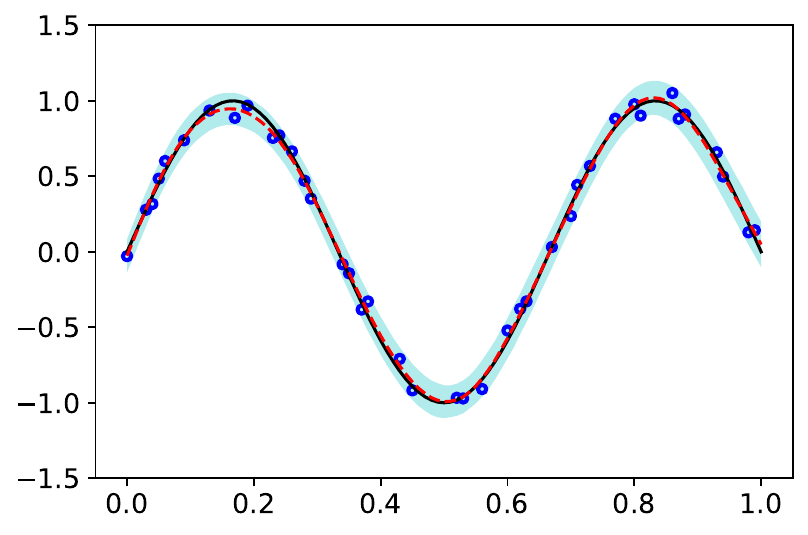}
			\caption{$f$}
		\end{subfigure}
		\hfill
		\begin{subfigure}[b]{0.3\textwidth}
			\centering
			\includegraphics[width=\textwidth, trim=0 0 0 0, clip]{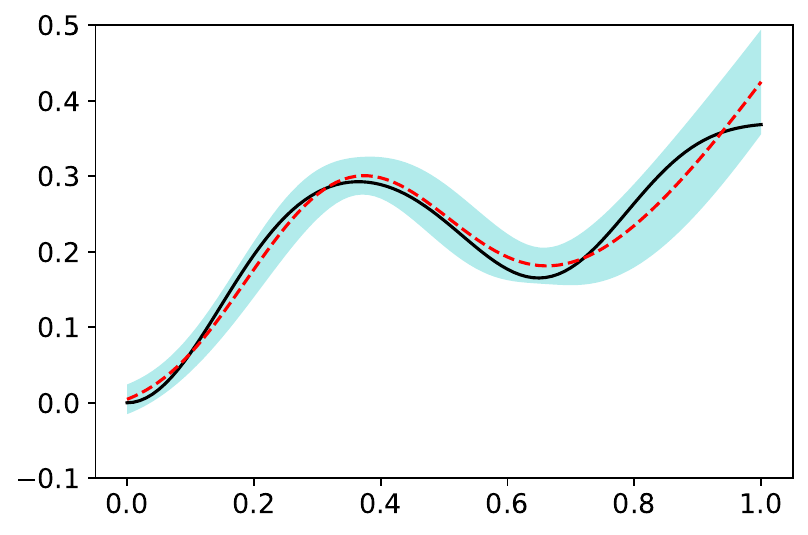}
			\caption{$\phi$}
		\end{subfigure}
		
		\caption{ODE system. Predicted $u$, $f$ and $\phi$ with model correction are based on a global integral operator structure. The training data, true solution, model mean prediction, and predictive uncertainty are respectively marked by blue circles, a solid black line, a red dashed line, and a shaded region spanning ±2 standard deviations.}
		\label{fig:IO_uands}
	\end{figure}

	We also consider another corrected architecture, in which the predicted mean solution and the correction term are generated by two separate integral-operator layers. The correction decoder follows the formulation in Eq.~\eqref{eq:IO_s}. The numerical results in Figure~\ref{fig:IO_uands} indicate that this variant exhibits behavior similar to that of the architecture presented above.
	
	\subsection{Reaction-Diffusion Equation}\label{Reaction_diffusion_equation}
	Next, we consider a time-dependent PDE that governs the dynamics of a reaction-diffusion system:
	\begin{equation}\label{eq_dr}
		\begin{aligned}
			\frac{\partial u}{\partial t} &= D \frac{\partial^{2} u}{\partial x ^{2}} + \lambda g(u) + f(x, t), \quad x \in [0, 1], \ t \in [0, 1], \\
			u(x, 0) &= 0.5 \sin^{2}(\pi x), \quad x \in [0, 1], \\
			u(0, t) &= u(1, t) = 0, \quad t \in [0, 1].
		\end{aligned}
	\end{equation}
	We set $D=0.01$, choose $g(u)=u(1-u)$, and use $\lambda=2$ in the ground-truth system. The source term is set to zero, i.e., $f(x,t)=0$. Following the same setting as in Section~\ref{ODE_System}, we consider three test cases:
	\begin{itemize}
		\item S1: The reaction is correctly specified as $g(u)=u(1-u)$, and $\lambda\ge 0$ is treated as unknown.
		\item S2: The reaction is misspecified as $g(u)=u^{2}$, and $\lambda\ge 0$ is treated as unknown.
		\item S3: The reaction is misspecified as $g(u)=u^{2}$ with $\lambda$ set to $1$. We augment the reaction term with a learnable correction, yielding $\lambda u^{2}+\mathcal{S}_{\psi}$.
	\end{itemize}
	For the three physical models in S1-S3, the estimated reaction terms are evaluated as
	$\tilde{\lambda}u_{\theta}(1-u_{\theta})$, $\tilde{\lambda}u_{\theta}^{2}$, and $\lambda u_{\theta}^{2} + \mathcal{S}_{\psi}$, respectively, where $\tilde{\lambda}$ is estimated by LVM-GP.
	In this example, we use $121$ measurements of $u$ and $195$ measurements of $f$, uniformly sampled over $(x,t)\in[0,1]\times[0,1]$.
	Both $u$ and $f$ are corrupted by additive Gaussian noise with standard deviations $0.02$ and $0.05$, respectively.
	The reference solution is generated by numerically solving Eq.~\eqref{eq_dr}.
	In S3, we apply the LVM-GP framework to infer $u$ and learn the correction term via $\mathcal{S}_{\psi}$.
	
	The encoder and solution decoder are implemented as fully connected networks with three hidden layers of width $128$ and Mish activations, while the correction decoder uses two hidden layers of width $64$ with Tanh activations.
	The regularization parameter $\beta$ is set to $0.5$. Training is conducted for $10,000$ iterations, divided into two phases: the first $5,000$ iterations focus solely on optimizing the predictive mean, while in the subsequent $5,000$ iterations, the predictive mean continues to be fine-tuned while the predictive standard deviation is simultaneously optimized. 
	
\begin{figure}[H]
		\begin{center}
			\begin{overpic}[width=0.22\textwidth, trim=0 0 0 0, clip=True]{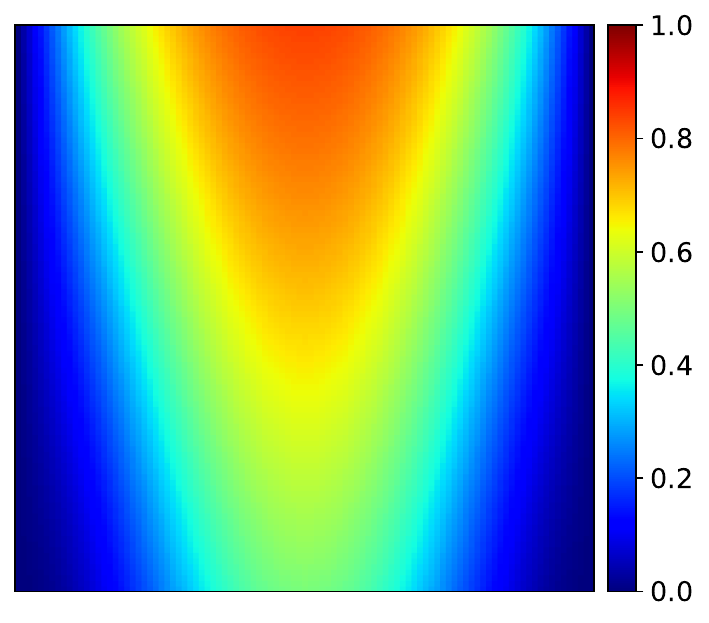}
				\put(35,90){\small{$u\big|$ref}}
			\end{overpic}
			\begin{overpic}[width=0.22\textwidth, trim=0 0 0 0, clip=True]{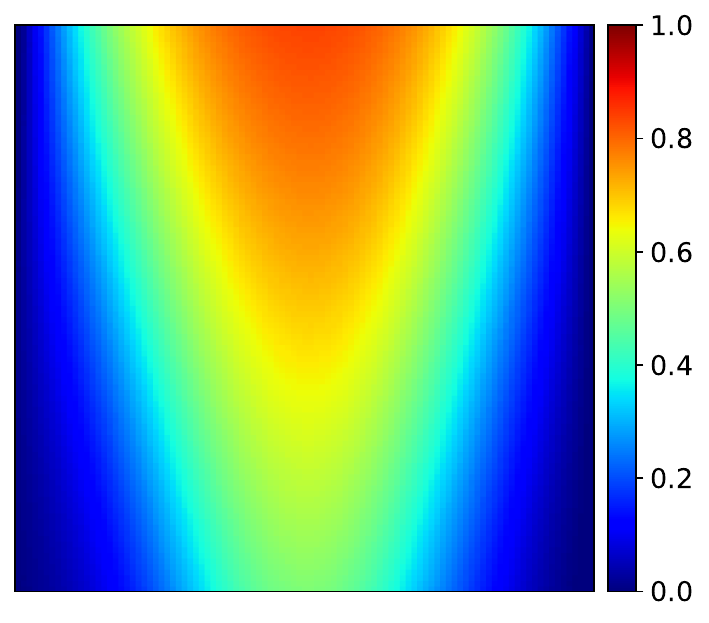}
				\put(20,90){\small{right model}}
			\end{overpic}
			\begin{overpic}[width=0.22\textwidth, trim=0 0 0 0, clip=True]{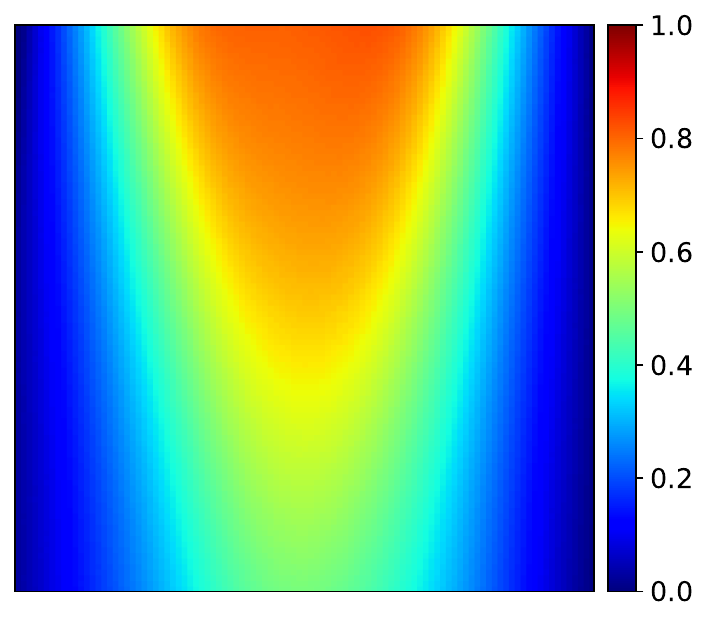}
				\put(5,90){\small{misspecified model}}
			\end{overpic}
			\begin{overpic}[width=0.22\textwidth, trim=0 0 0 0, clip=True]{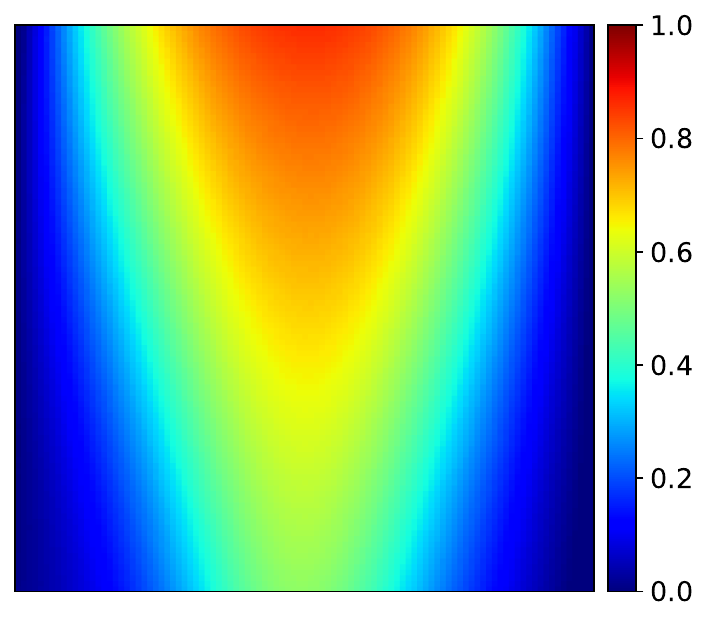}
				\put(10,90){\small{corrected model}}
			\end{overpic}
			\subcaption{Predictions for $u$(mean).}
			
			\vspace{20pt}
			
			\begin{overpic}[width=0.22\textwidth, trim=0 0 0 0, clip=True]{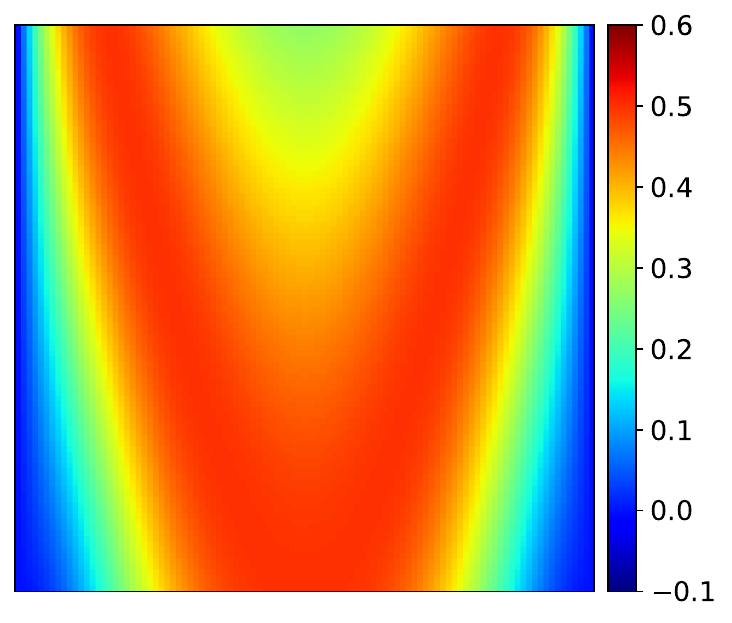}
				\put(35,87.5){\small{$\phi\big|$ref}}
			\end{overpic}
			\begin{overpic}[width=0.22\textwidth, trim=0 0 0 0, clip=True]{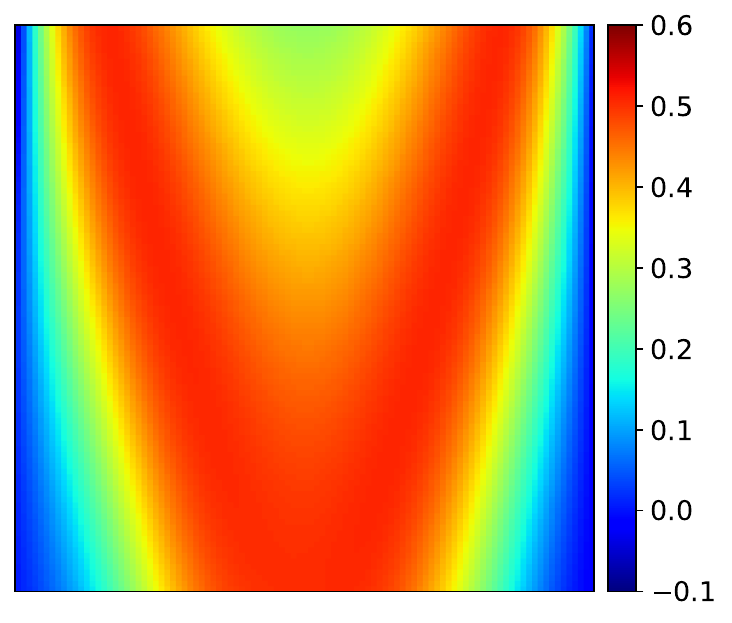}
				\put(20,87.5){\small{right model}}
			\end{overpic}
			\begin{overpic}[width=0.22\textwidth, trim=0 0 0 0, clip=True]{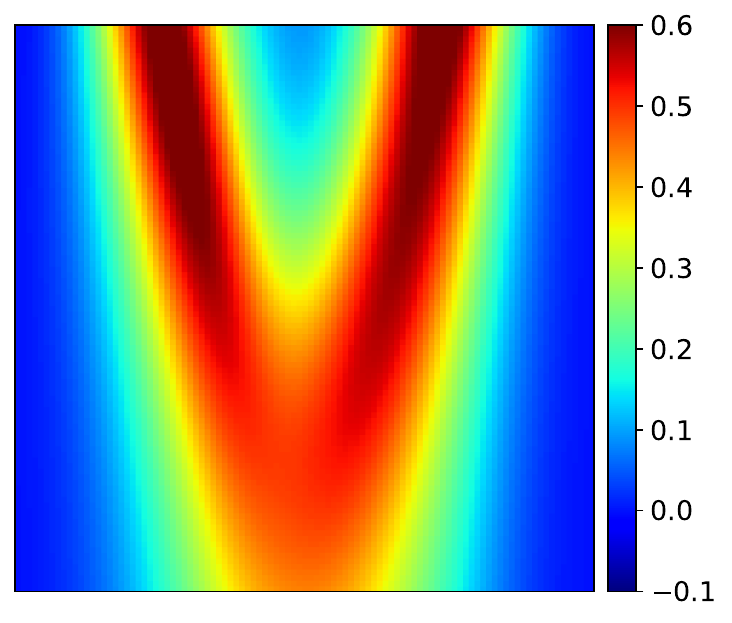}
				\put(5,87.5){\small{misspecified model}}
			\end{overpic}
			\begin{overpic}[width=0.22\textwidth, trim=0 0 0 0, clip=True]{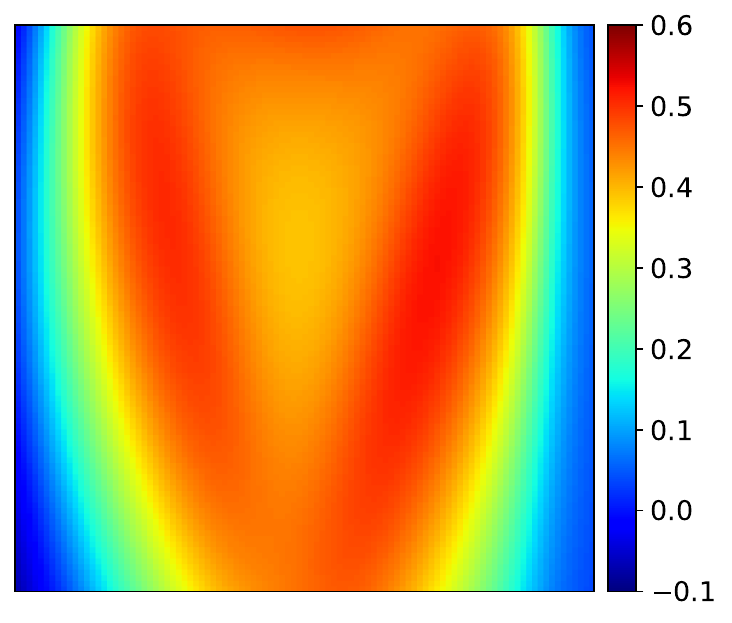}
				\put(10,87.5){\small{corrected model}}
			\end{overpic}
			\subcaption{Predictions for $\phi$(mean).}
		\end{center}
		\caption{Time-dependent reaction-diffusion system. The predicted means of $u$ and $\phi$ are obtained using the LVM-GP (S1 and S2) and the LVM-GP-based model correction framework (S3), applied to noisy observations of $u$ and $f$ under different reaction models. See Figure~\ref{fig:dr_cross_section} for uncertainty quantification.}
		\label{fig:dr}
	\end{figure}

\begin{figure}[H]
		\begin{center}
			\begin{overpic}[width=0.3\textwidth, trim=0 0 0 0, clip=True]{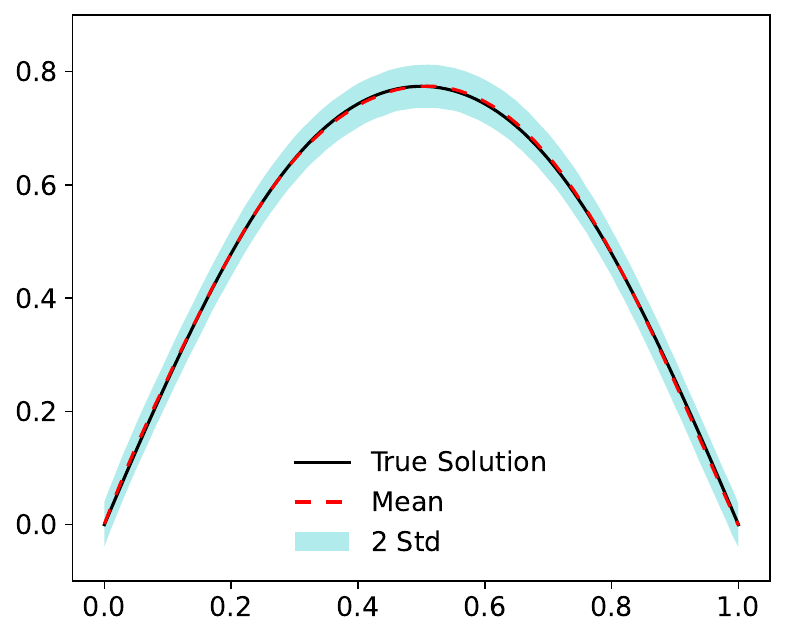}
				\put(52.5,-2.5){\small{$x$}}
				\put(-5, 41.5) {\small{$u$}}
				\put(35,82.5){\small{right model}}
			\end{overpic}
			\hspace{5pt}
			\begin{overpic}[width=0.3\textwidth, trim=0 0 0 0, clip=True]{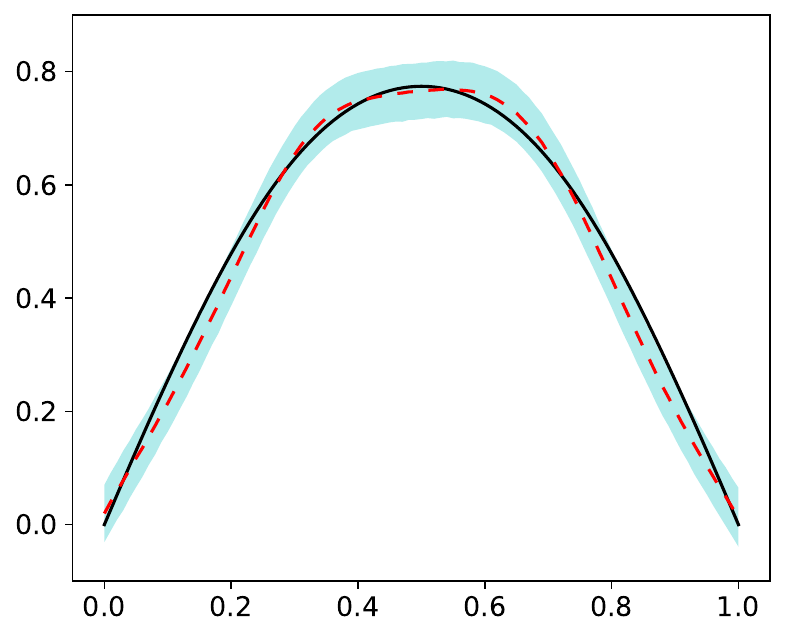}
				\put(52.5,-2.5){\small{$x$}}
				\put(23.5,82.5){\small{misspecified model}}
			\end{overpic}
			\hspace{5pt}
			\begin{overpic}[width=0.3\textwidth, trim=0 0 0 0, clip=True]{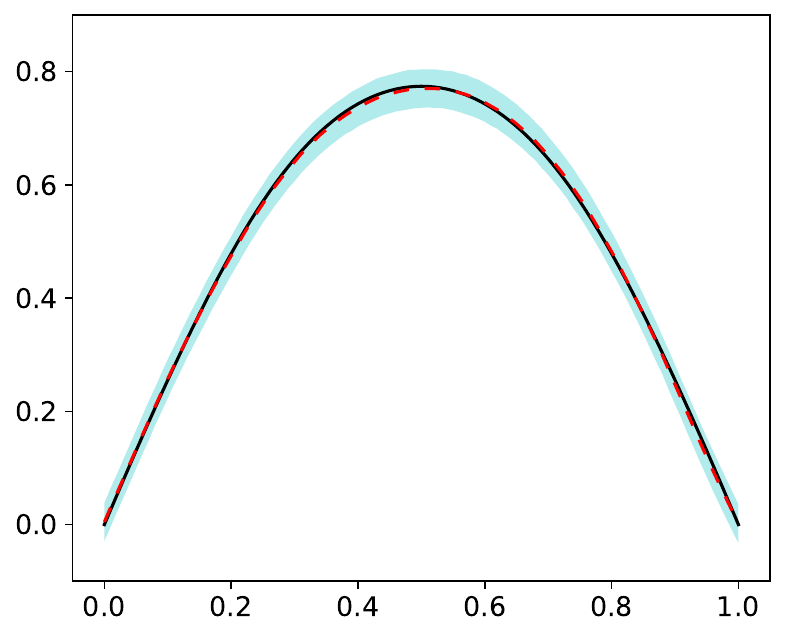}
				\put(52.5,-2.5){\small{$x$}}
				\put(28.5,82.5){\small{corrected model}}
			\end{overpic}
			\subcaption{Predictions for $u(t = 0.75)$.}
			
			\vspace{20pt}
			
			\begin{overpic}[width=0.305\textwidth, trim=0 0 0 0, clip=True]{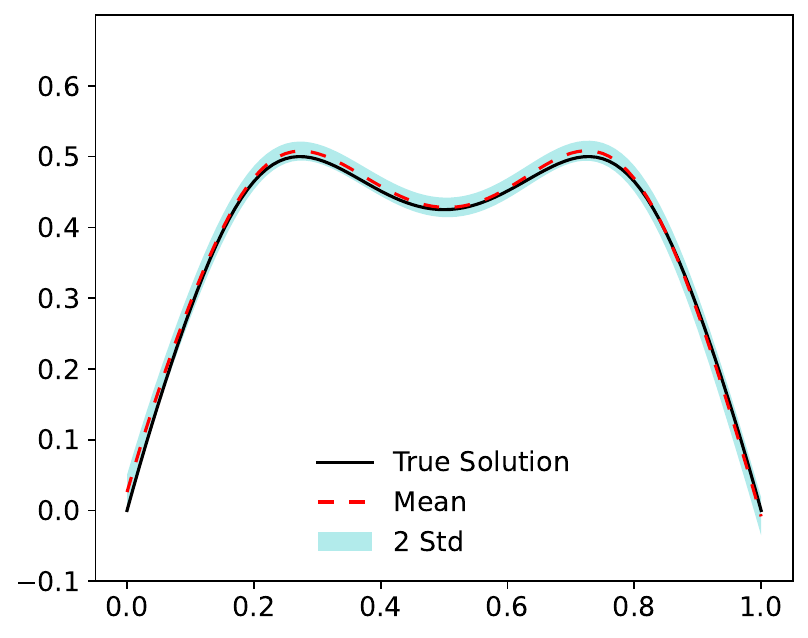}
				\put(53.5,-2.5){\small{$x$}}
				\put(-3, 42.5) {\small{$\phi$}}
				\put(37.5,80){\small{right model}}
			\end{overpic}
			\hspace{5pt}
			\begin{overpic}[width=0.305\textwidth, trim=0 0 0 0, clip=True]{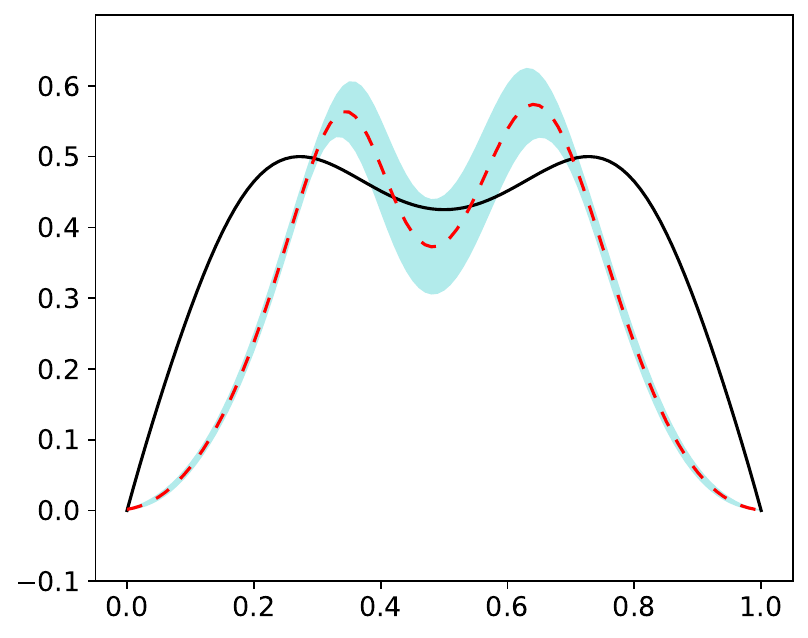}
				\put(53.5,-2.5){\small{$x$}}
				\put(25,80){\small{misspecified model}}
			\end{overpic}
			\hspace{3pt}
			\begin{overpic}[width=0.305\textwidth, trim=0 0 0 0, clip=True]{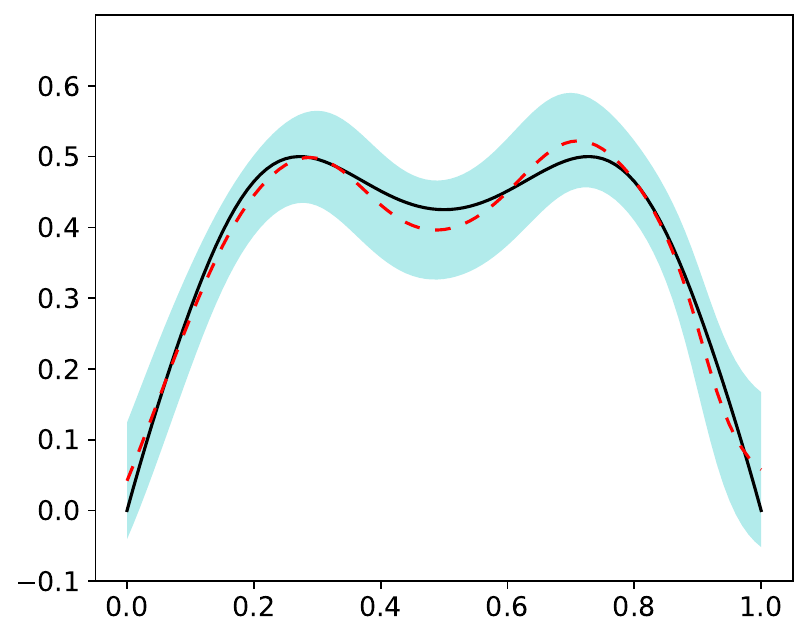}
				\put(53.5,-2.5){\small{$x$}}
				\put(30,80){\small{corrected model}}
			\end{overpic}
			\subcaption{Predictions for $\phi(t = 0.5)$.}
		\end{center}
		\caption{Time-dependent reaction-diffusion system in slices $t = 0.75$ for $u$, $t = 0.5$ for $\phi$. The true solution, model mean prediction, and predictive uncertainty are respectively marked by a solid black line, a red dashed line, and a shaded region spanning ±2 standard deviations.}
		\label{fig:dr_cross_section}
	\end{figure}
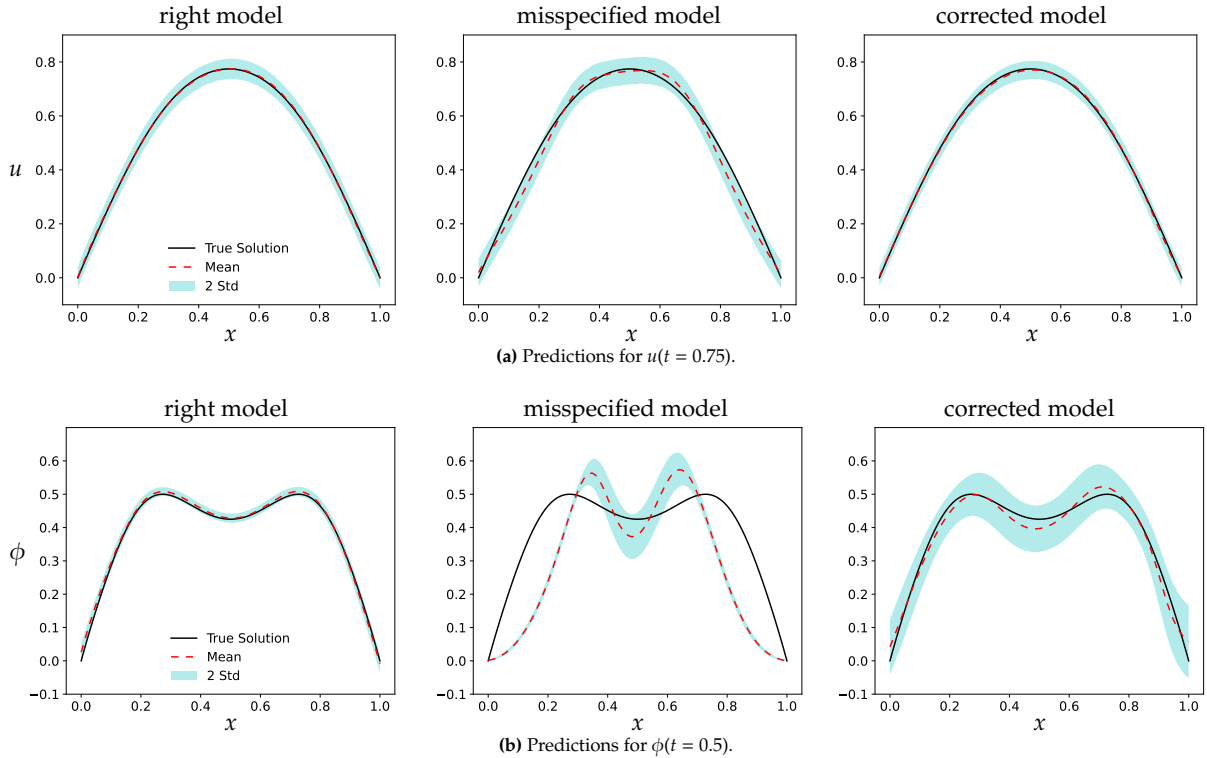

	As shown in Figure~\ref{fig:dr} and Figure~\ref{fig:dr_cross_section}, when the model is correctly specified (S1), the predicted $u$ and $\phi$ match well with the reference solutions, and all errors lie within the uncertainty bounds, indicating the reliability of the results. When the model is misspecified (S2), the predicted results deviate significantly from the references, especially for the reaction term. However, after correcting the misspecified model using our proposed method (S3), the predictions improve substantially compared to S2. Although the accuracy remains lower than in S1, the prediction errors consistently fall within the uncertainty region, demonstrating that the results are still reliable. It is worth noting that the uncertainty in S3 is noticeably larger than that in S1, reflecting the additional model uncertainty introduced by the correction network.
	
	\subsection{2D Non-Newtonian Flows}
	In this example, we evaluate the proposed approach for non-Newtonian flows. The method starts with a misspecified Newtonian constitutive relation, which is then corrected by a neural network. We investigate two steady-state configurations: a 2D channel flow and a 2D cavity flow, both using power-law fluids. The constitutive relation for a power-law fluid is $\mu = \mu_{0}|\bm{S}|^{n-1}$. Here, $\mu$ is the dynamic viscosity, $\mu_{0}$ is a constant, $\bm{S} = \nabla \bm{u} + (\nabla \bm{u})^{T}$ is the shear rate tensor (with $|\bm{S}| = \sqrt{(\bm{S}:\bm{S})/2}$), $\bm{u}$ is the velocity, and $n$ is the power-law index. Specifically, shear-thinning (pseudoplastic, $n < 1$) and shear-thickening (dilatant, $n > 1$) fluids are considered for the channel and cavity flows, respectively.
	\subsubsection{Channel Flow}\label{channel_flow}
	We first consider a 2D non-Newtonian channel flow studied in \cite{wang2015localized}:
	\begin{equation}\label{eq_channels}
		\begin{aligned}
			&\frac{\partial }{\partial y}\left(\mu(y)\frac{\partial u}{\partial y}\right) - \frac{\partial P}{\partial x} = f(y), \quad y \in \left[-\frac{H}{2}, \frac{H}{2}\right], \\
			&u\left(-\frac{H}{2}\right) = u\left(\frac{H}{2}\right) = 0,
		\end{aligned}
	\end{equation}
	with $H=1$, a constant pressure gradient $\partial P/ \partial x=c$, and $f(y)=0$.
	Under these conditions, the analytical solution is
	\begin{equation*}
		u(y) = \frac{n}{n+1}\left(-\frac{1}{\mu_{0}}\frac{\partial P}{\partial x}\right)^{1/n}\left[\left(\frac{H}{2}\right)^{1+1/n} - |y|^{1+1/n}\right].
	\end{equation*}
	For generating the reference solution and training data, we set $n=0.25$, $c=-1$, and $\mu_{0}=0.5$. We consider two data regimes with identical sample sizes: a noise-free setting and a noisy setting. In both regimes, we randomly and uniformly sample $30$ measurements of $u$ and $51$ measurements of $f$ from $[-H/2,H/2]$.
	
	In many applications, the true fluid is non-Newtonian, but the exact constitutive relation is not available. It is therefore common to approximate the system by a Newtonian model with a spatially uniform viscosity. Following this idea, we first consider the Newtonian model, where the viscosity is an unknown constant $\mu_1$:
	\begin{equation}\label{eq:eq_channels_flow_Newtonian_flows}
		\frac{\partial}{\partial y}\left(\mu_1 \frac{\partial u}{\partial y}\right) - \frac{\partial P}{\partial x} = f(y).
	\end{equation}
	Given the measurements, we apply LVM-GP to jointly infer the constant viscosity $\mu_1$ and the velocity profile $u(y)$ under this Newtonian model.
	
	Even after estimating $\mu_1$, a constant-viscosity Newtonian model generally cannot reproduce the behavior of the non-Newtonian model truth, so a systematic model discrepancy remains. To represent and correct this discrepancy, we next consider the corrected Newtonian model, which adds an unknown correction function $s(y)$ to the Newtonian equation:
	\begin{equation}\label{eq_channels_flow_corrected}
		\frac{\partial }{\partial y}\left(\mu_{1}\frac{\partial u}{\partial y}\right) - \frac{\partial P}{\partial x} + s(y) = f(y).
	\end{equation}
	In this corrected formulation, we fix the viscosity to $\mu_1=0.1$ and apply LVM-GP to jointly learn the correction term $s(y)$ together with the solution $u(y)$. In short, the first model accounts for mismatch only through the single constant parameter $\mu_1$, whereas the corrected model keeps $\mu_1$ fixed and uses the learned function $s(y)$ to compensate for the missing non-Newtonian effects.
	
	For all experiments in this subsection, the encoder and solution decoder are implemented as three-layer fully connected networks of width $64$ with Mish activations. In the noise-free setting, both the correction network and the constitutive-relation DNN are two-hidden-layer MLPs of width $32$ with Tanh activations, with $\beta=10^{-6}$. In the noisy setting, the correction network uses a two-hidden-layer MLP of width $50$ with Tanh activations and $\beta=0.5$. The training schedule and number of iterations follow Section~\ref{ODE_System}.
	
\begin{figure}[h]
		\begin{center}
			\begin{overpic}[width=0.33\textwidth, trim=0 0 0 0, clip=True]{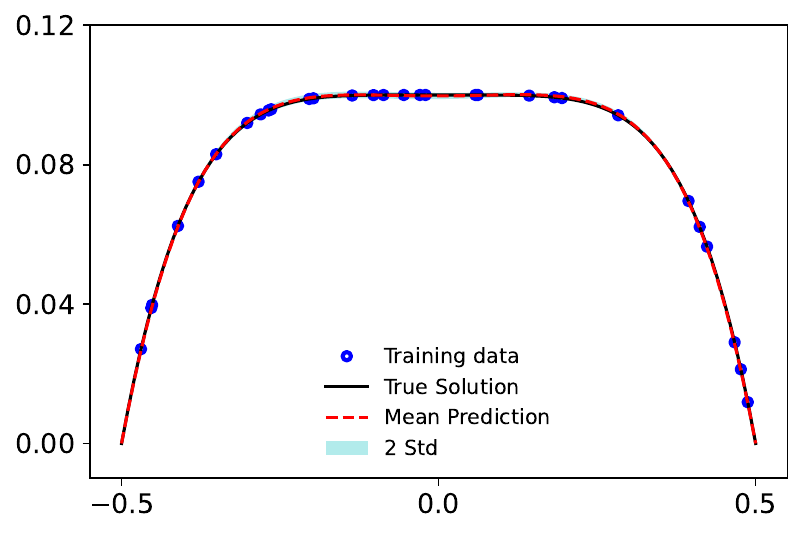}
				\put(52.5, -2.5){\small{$y$}}
				\put(-5, 32.5) {\small{$u$}}
			\end{overpic}    
			\hspace{10pt}
			\begin{overpic}[width=0.33\textwidth, trim=0 0 0 0, clip=True]{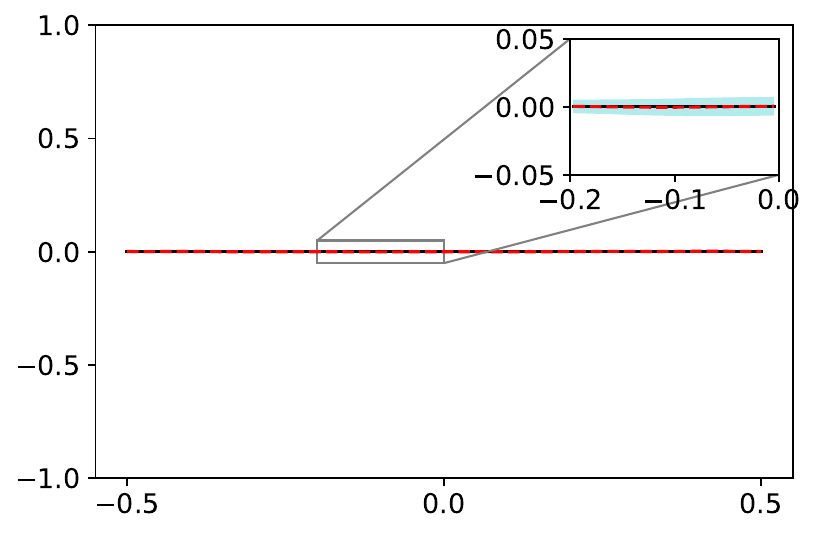}
				\put(52.5, -2.5){\small{$y$}}
				\put(-1.5, 32.5) {\small{$f$}}
			\end{overpic}
			\subcaption{$\mu_{1}\frac{\partial^{2}{u}}{\partial y^{2}} - \frac{\partial P}{\partial x} + s(y) = f(y)$}
			
			\vspace{20pt}
			
			\begin{overpic}[width=0.33\textwidth, trim=0 0 0 0, clip=True]{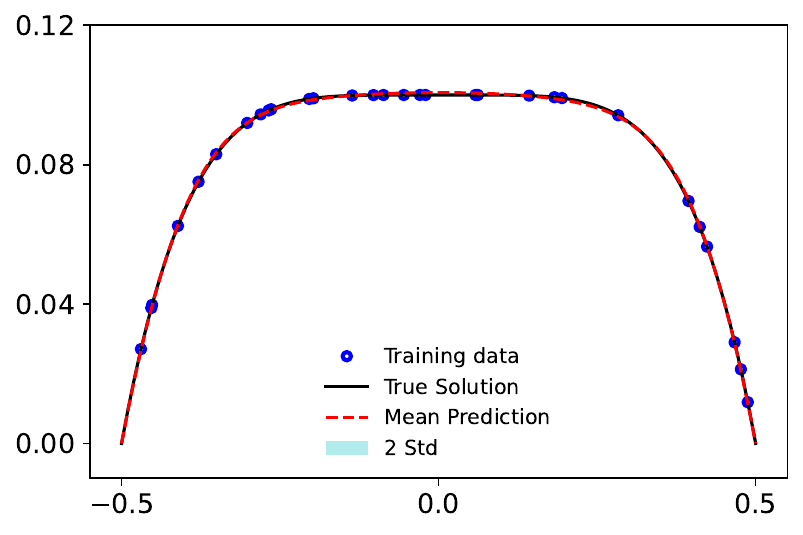}
				\put(52.5, -2.5){\small{$y$}}
				\put(-5, 32.5) {\small{$u$}}
			\end{overpic}
			\hspace{10pt}
			\begin{overpic}[width=0.33\textwidth, trim=0 0 0 0, clip=True]{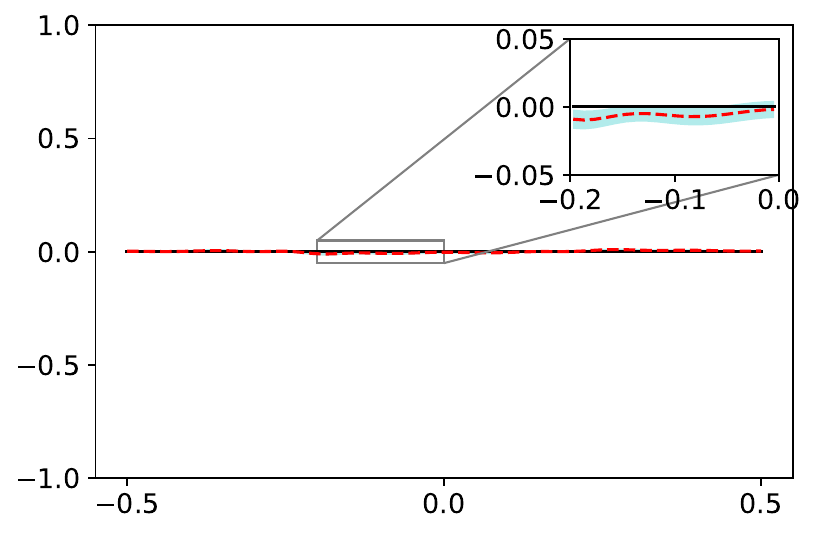}
				\put(52.5, -2.5){\small{$y$}}
				\put(-1.5, 32.5) {\small{$f$}}
			\end{overpic}
			\subcaption{$\frac{\partial}{\partial y}(\mu(y)\frac{\partial u}{\partial y}) - \frac{\partial P}{\partial x} = f(y)$}
		\end{center}
		\caption{2D non-Newtonian channel flow. Predictive results for $u$ and $f$. A comparison of two approaches for correcting misspecified physics: (a) modeling the equation discrepancy $s$ via a correction network, and (b) learning the spatially varying viscosity $\mu(y)$ using a deep neural network (DNN). The training data, true solution, model mean prediction, and predictive uncertainty are respectively marked by blue circles, a solid black line, a red dashed line, and a shaded region spanning ±2 standard deviations.}
		\label{fig:channels_flow_noisefree}
	\end{figure}
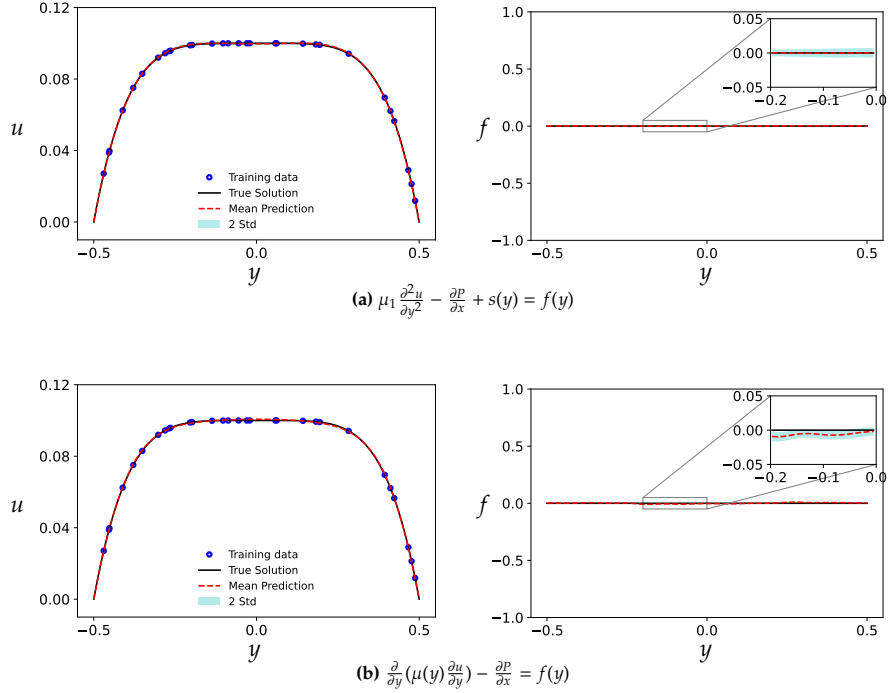

	We first consider the noise-free regime and compare two strategies for handling the misspecification.
	Our approach learns the discrepancy term $s(y)$ in Eq.~\eqref{eq_channels_flow_corrected}, whereas the baseline directly learns the constitutive relation by treating the viscosity in Eq.~\eqref{eq:eq_channels_flow_Newtonian_flows} as an unknown function $\mu(y)$ parameterized by a DNN.
	As shown in Figure~\ref{fig:channels_flow_noisefree}, both methods accurately predict $u$, but our approach yields a noticeably smoother and more accurate prediction for $f$.
	
	We next turn to the noisy regime. The sample sizes for $u$ and $f$ are unchanged, but the data are contaminated by additive Gaussian noise with standard deviations $0.005$ and $0.05$, respectively.
	To assess the benefit of the correction network, we compare (i) fitting the misspecified Newtonian model Eq.~\eqref{eq:eq_channels_flow_Newtonian_flows} and inferring the constant viscosity $\mu_{1}$ from noisy data, and (ii) fixing $\mu_{1}=0.1$ and learning $s(y)$ in the corrected model Eq.~\eqref{eq_channels_flow_corrected}.
	
\begin{figure}[H]
		\begin{center}
			\begin{overpic}[width=0.3\textwidth, trim=0 0 0 0, clip=True]{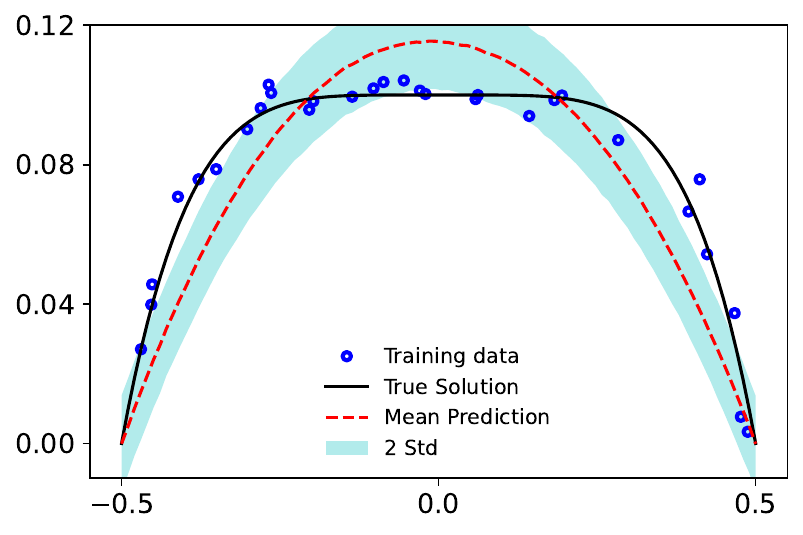}
				\put(52.5, -2.5){\small{$y$}}
				\put(-5, 32.5) {\small{$u$}}
			\end{overpic}
			\hspace{10pt}
			\begin{overpic}[width=0.3\textwidth, trim=0 0 0 0, clip=True]{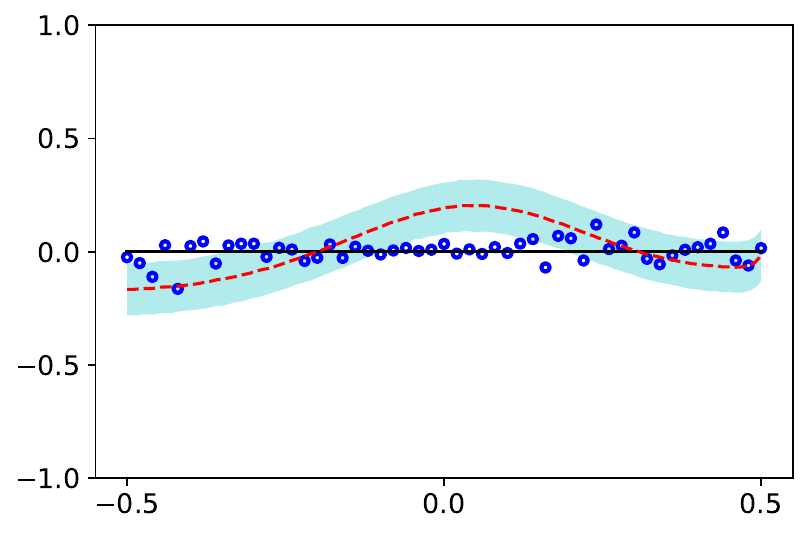}
				\put(52.5, -2.5){\small{$y$}}
				\put(-2.5, 32.5) {\small{$f$}}
			\end{overpic}
			\vspace{0.5em}
			\subcaption{LVM-GP with a misspecified physical model}
			
			\vspace{10pt}
			
			\begin{overpic}[width=0.3\textwidth, trim=0 0 0 0, clip=True]{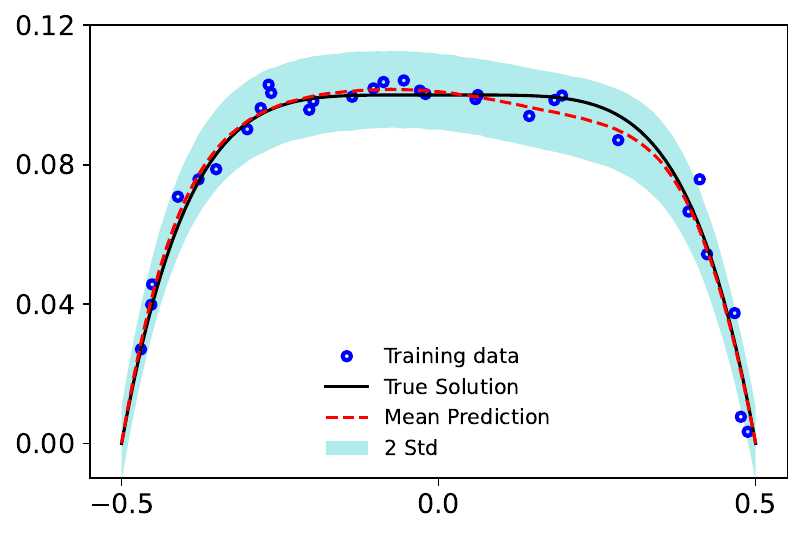}
				\put(52.5, -2.5){\small{$y$}}
				\put(-5, 32.5) {\small{$u$}}
			\end{overpic}
			\hspace{10pt}
			\begin{overpic}[width=0.3\textwidth, trim=0 0 0 0, clip=True]{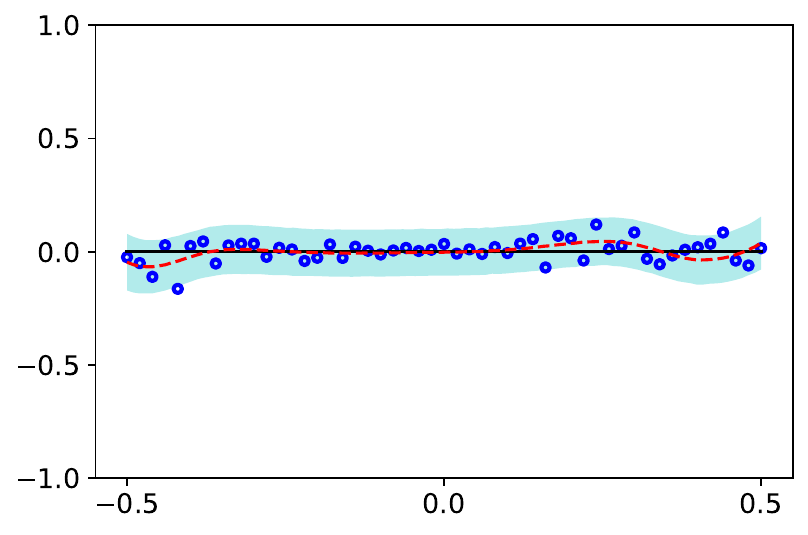}
				\put(52.5, -2.5){\small{$y$}}
				\put(-2.5, 32.5) {\small{$f$}}
			\end{overpic}
			\hspace{10pt}
			\begin{overpic}[width=0.3\textwidth, trim=0 0 0 0, clip=True]{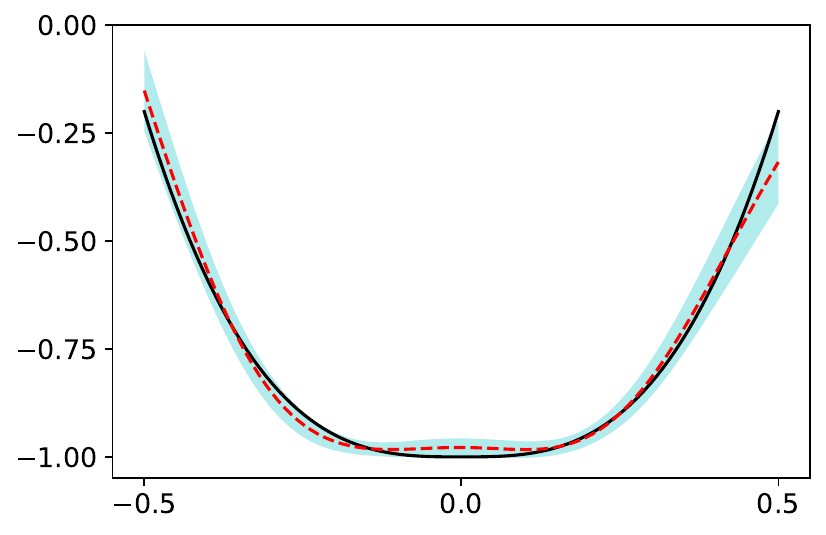}
				\put(52.5, -2.5){\small{$y$}}
				\put(-2.5, 32.5) {\small{$s$}}
			\end{overpic}
			\vspace{0.5em}
			
			\subcaption{Correcting the misspecified physical model by a correction network using LVM-GP}
		\end{center}
		\caption{2D non-Newtonian channel flow. Predictive results for $u$, $f$, and the correction term $s$ with noisy measurement data under two modeling scenarios: (a) LVM-GP with a misspecified physical model and (b) correcting the misspecified physical model by a correction network using LVM-GP. The training data, true solution, model mean prediction, and predictive uncertainty are respectively marked by blue circles, a solid black line, a red dashed line, and a shaded region spanning ±2 standard deviations.}
		\label{fig:channels_flow_noise}
	\end{figure}
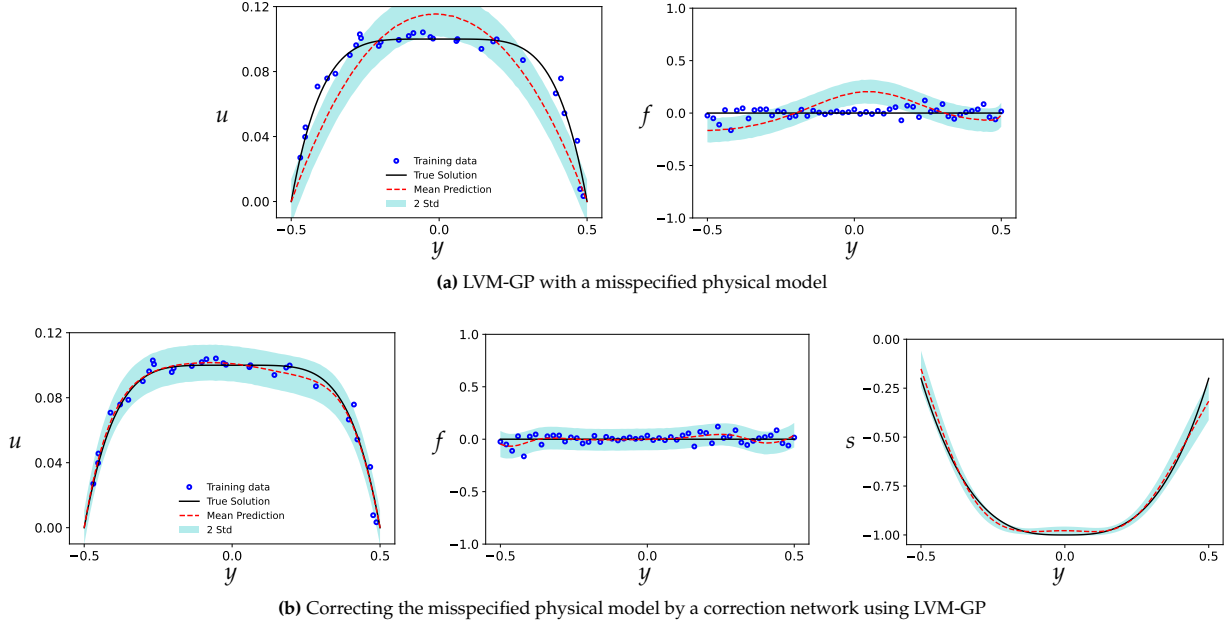

	For this flow, the correction term admits a closed-form expression. Substituting the power-law solution into Eq.~\eqref{eq_channels_flow_corrected} yields
	$s(y)=4\mu_{1}c^{4}|y|^{3}/\mu^{4}_{0}+c$,
	which serves as a reference for evaluating the learned correction.
	Figure~\ref{fig:channels_flow_noise} shows that the correction module improves performance substantially relative to the misspecified Newtonian fit. Moreover, the prediction errors for both the solution and the discrepancy term are covered by the predictive uncertainty bounds, indicating reliable uncertainty quantification.
	
	\subsubsection{Cavity flow}
	We next consider steady non-Newtonian flow in a two-dimensional lid-driven square cavity ($x,y\in[0,1]$), governed by the incompressible Navier--Stokes equations
	\begin{subequations}
		\begin{align}
			\nabla \cdot \bm{u} &= 0, \label{eq:continuity}\\
			\bm{u}\cdot \nabla \bm{u} &= -\nabla P + \nabla \cdot \big[\nu(\nabla \bm{u}+(\nabla \bm{u})^{T})\big] + \bm{f}. \label{eq:momentum_eq}
		\end{align}
	\end{subequations}
	Here, $\bm{u}=(u,v)$ denotes the velocity field, $P$ the pressure, and $\nu$ the kinematic viscosity. On the upper lid, we prescribe
	\[
	u = U\Big(1-\frac{\cosh[r(x-\frac{L}{2})]}{\cosh(\frac{rL}{2})}\Big), \qquad v=0,
	\]
	with $U=0.01$, $r=10$, and $L=1$, while the remaining walls satisfy no-slip conditions. The viscosity follows the same power-law model as in the channel-flow example, $\mu=\rho\nu=\mu_0|\bm{S}|^{n-1}$. Training data and reference solutions are generated using the lattice Boltzmann equation (LBE) solver from \cite{wang2015localized} with power-law index $n=1.5$. We follow \cite{zou2024correcting} to generate the non-Newtonian cavity-flow dataset.
	
	As in Section~\ref{channel_flow}, an incorrect constitutive relation can be compensated by introducing a correction source term:
	\begin{subequations}
		\begin{align}
			\nabla \cdot \bm{u} &= 0, \label{eq:correct_continuous}\\
			\bm{u}\cdot \nabla \bm{u} &= -\nabla P + \nu_{0}\nabla^{2}\bm{u} + \bm{f} + \bm{s}, \label{eq:correct_Rom}
		\end{align}
	\end{subequations}
	where $\bm{s}=(s_x,s_y)$ corrects the $x$- and $y$-momentum equations, respectively, and $\bm{f}\equiv \bm{0}$ on $[0,1]^2$.
	
	For this example, we use 400 randomly sampled measurements for each velocity component ($u$ and $v$) and pressure measurements on a uniform $81\times81$ grid. We compare three settings: (1) a baseline misspecified model, where the viscosity in Eq.~\eqref{eq:correct_Rom} is treated as an unknown constant and inferred by LVM-GP without correction; (2) Ensemble PINNs \cite{zou2024correcting}, where the viscosity is fixed at $\nu_0=10^{-4}$; and (3) our model-corrected LVM-GP, which also fixes $\nu_0=10^{-4}$ but learns $\bm{s}$ via our correction network.
	
	The encoder, solution decoder, and correction decoder are implemented as three-layer MLPs with 128 hidden units per layer. The encoder and solution decoder use the Mish activation, whereas the correction decoder uses Tanh. We train the model with Adam for 50,000 iterations, starting from a learning rate of $10^{-3}$ and decaying it by a factor of 0.9 every 1,000 steps. The hyperparameter $\beta$ is set to $10^{-8}$.

\begin{figure}[H]
		\begin{center}
			\begin{overpic}[width=0.32\textwidth, trim=0 0 0 0, clip=True]{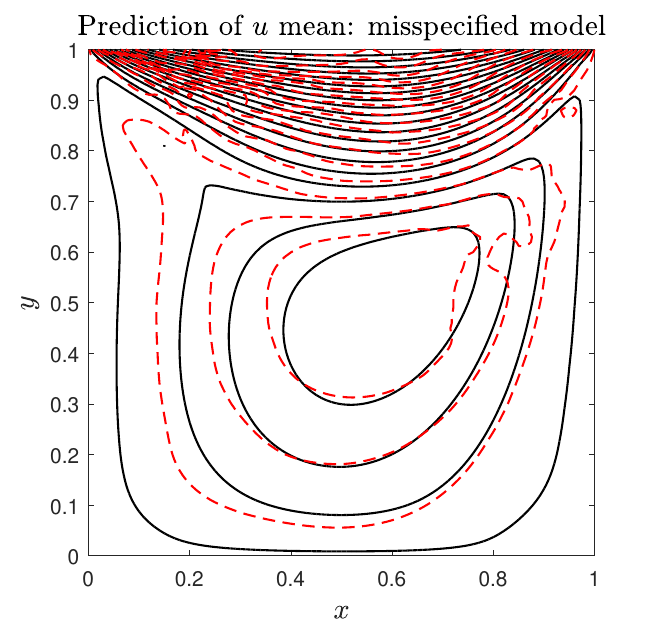}
			\end{overpic}
			\hspace{5pt}
			\begin{overpic}[width=0.32\textwidth, trim=0 0 0 0, clip=True]{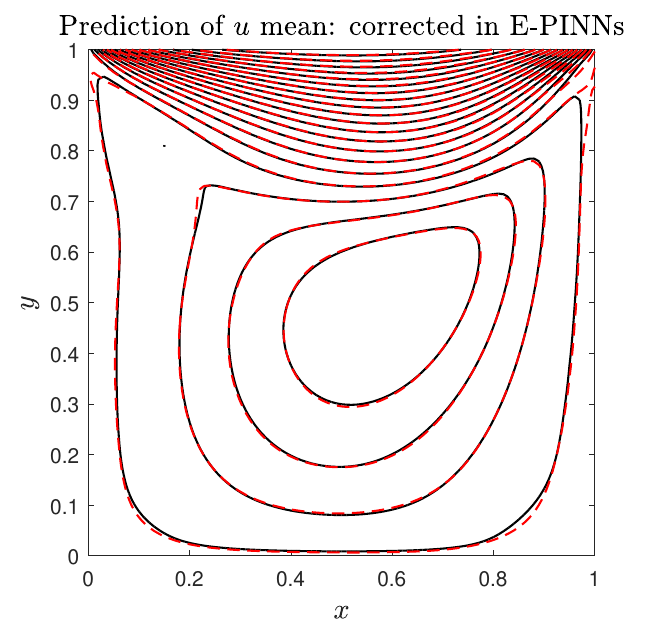}
			\end{overpic}
			\hspace{5pt}
			\begin{overpic}[width=0.32\textwidth, trim=0 0 0 0, clip=True]{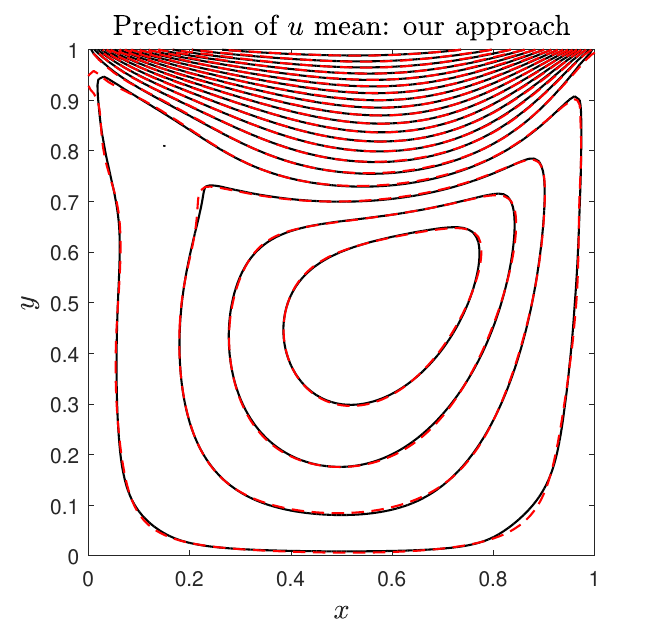}
			\end{overpic}
			
			\vspace{10pt}
			
			\begin{overpic}[width=0.32\textwidth, trim=0 0 0 0, clip=True]{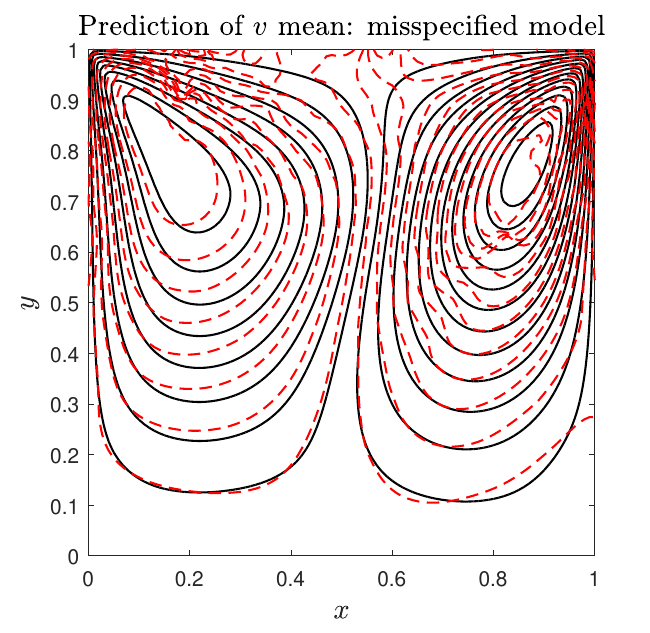}
				\put(5,-5){{ \footnotesize  (a) LVM-GP with misspecified model}}
			\end{overpic}
			\hspace{5pt}
			\begin{overpic}[width=0.32\textwidth, trim=0 0 0 0, clip=True]{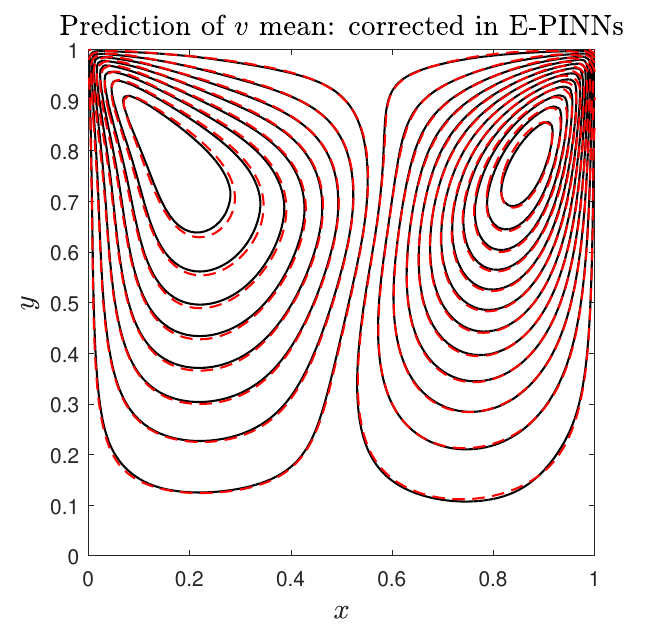}
				\put(5,-5){{ \footnotesize (b) Ensemble PINNs with model correction}}
			\end{overpic}
			\hspace{3pt}
			\begin{overpic}[width=0.32\textwidth, trim=0 0 0 0, clip=True]{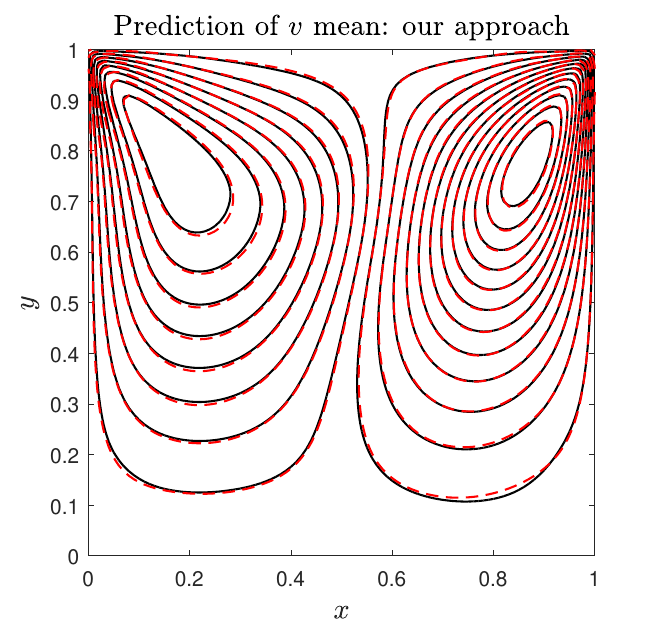}
				\put(10,-5){{ \footnotesize (c) LVM-GP with model correction}}
			\end{overpic}
		\end{center}
		\caption{Non-Newtonian cavity flow. Contour predictions of $u$ and $v$ for LVM-GP with model misspecification, Ensemble PINNs with model correction, and LVM-GP with model correction. The red dashed lines represent the predicted mean, while the black solid lines denote the reference solutions.} 
		\label{fig:cavity flow uv}
	\end{figure}

\begin{figure}[H]
		\begin{center}
			\begin{overpic}[width=0.32\textwidth, trim=0 0 0 0, clip=True]{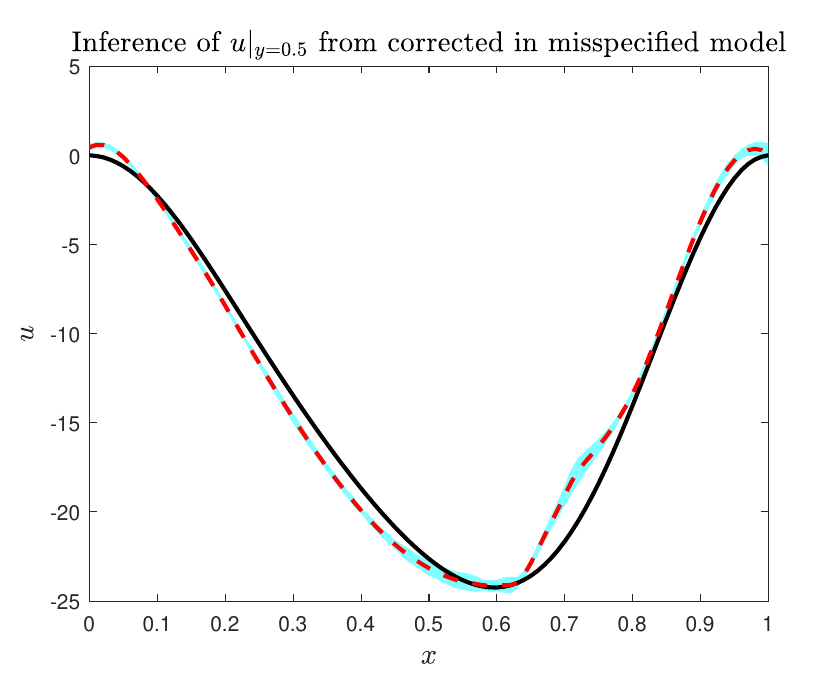}
			\end{overpic}
			\hspace{5pt}
			\begin{overpic}[width=0.32\textwidth, trim=0 0 0 0, clip=True]{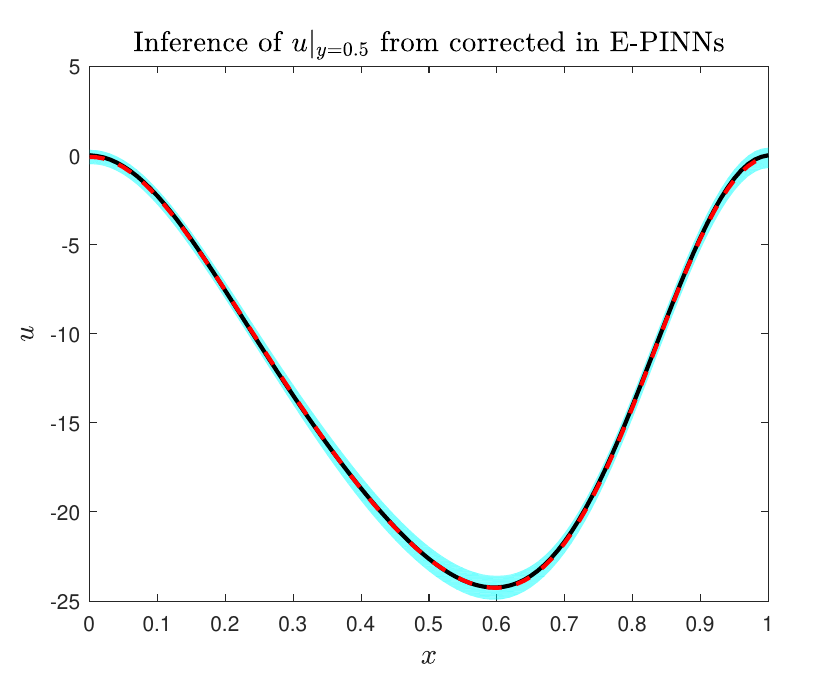}
			\end{overpic}
			\hspace{5pt}
			\begin{overpic}[width=0.32\textwidth, trim=0 0 0 0, clip=True]{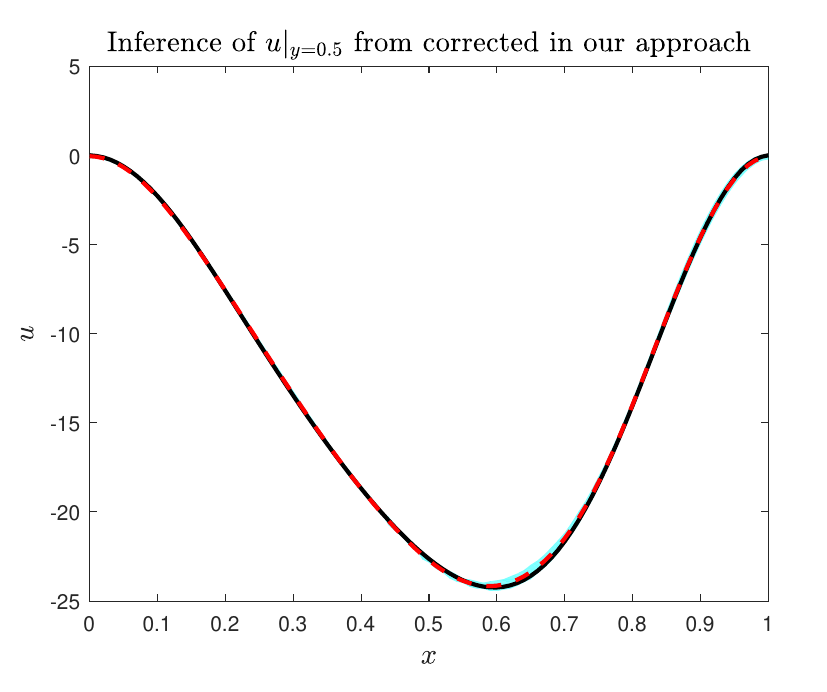}
			\end{overpic}
			
			\vspace{10pt}
			
			\begin{overpic}[width=0.32\textwidth, trim=0 0 0 0, clip=True]{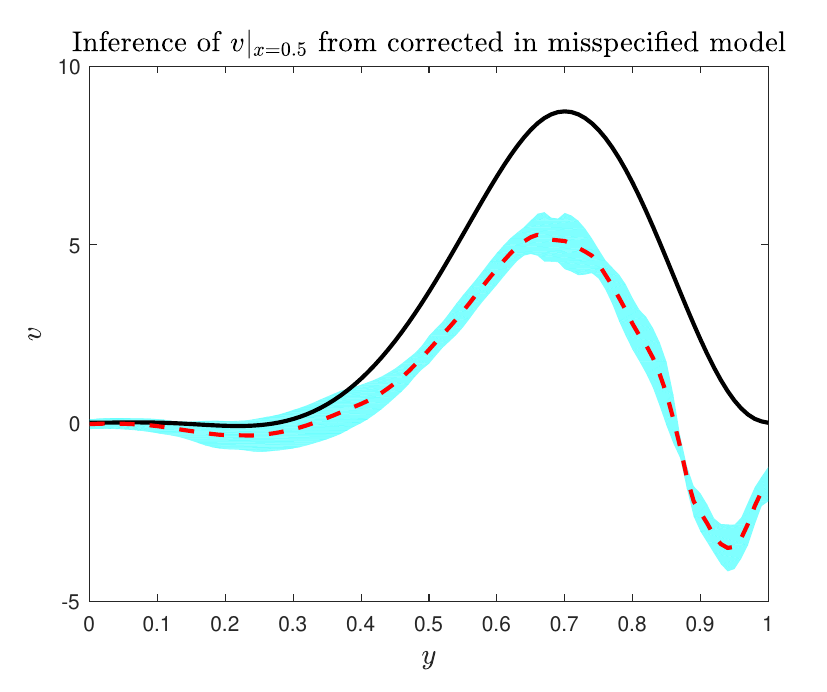}
				\put(5,-5){{ \footnotesize  (a) LVM-GP with misspecified model}}
			\end{overpic}
			\hspace{5pt}
			\begin{overpic}[width=0.32\textwidth, trim=0 0 0 0, clip=True]{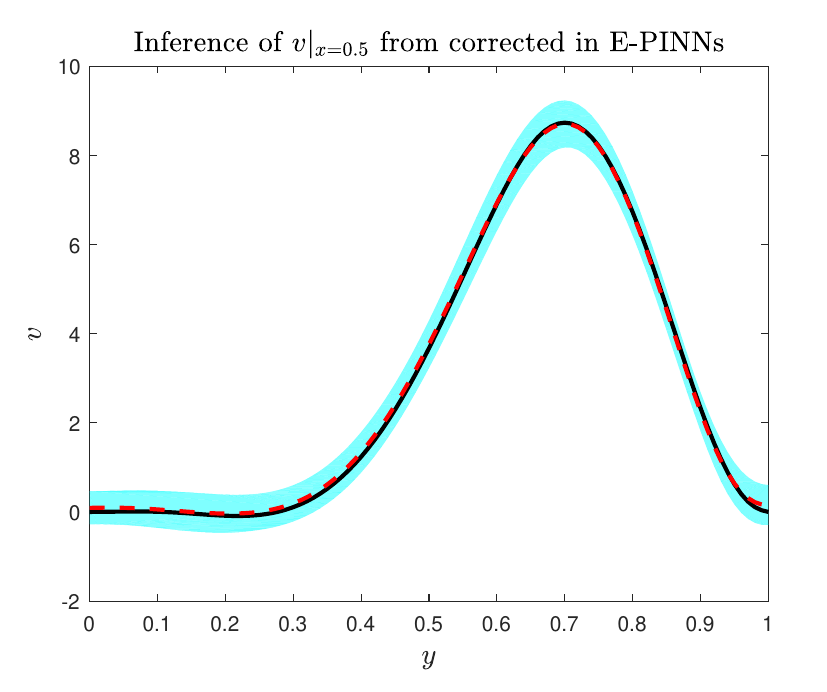}
				\put(5,-5){{ \footnotesize (b) Ensemble PINNs with model correction}}
			\end{overpic}
			\hspace{3pt}
			\begin{overpic}[width=0.32\textwidth, trim=0 0 0 0, clip=True]{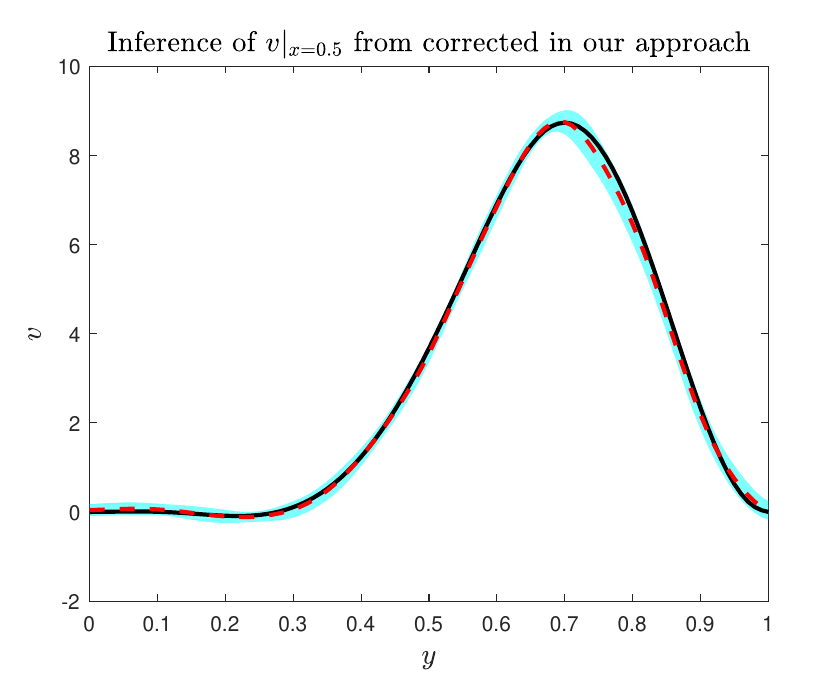}
				\put(10,-5){{ \footnotesize (c) LVM-GP with model correction}}
			\end{overpic}
		\end{center}
		\caption{Non-Newtonian cavity flow. 1D slice predictions of $u$ and $v$. The red dashed line represents the predicted mean, the black solid line denotes the reference solution, and the blue band indicates the uncertainty (twice the standard deviation).} 
		\label{fig:cavity flow uv slices}
	\end{figure}

	Figure~\ref{fig:cavity flow uv} compares the predicted $u$ and $v$ fields obtained by different methods. The misspecified LVM-GP model yields noticeably less accurate predictions, whereas the corrected model produces results that are much closer to the reference solution and comparable to those of Ensemble PINNs. This behavior is also reflected in the relative $L_2$ errors reported in Table~\ref{tab:cavity_flow}. Similar observations can be made from the centerline profiles in Figure~\ref{fig:cavity flow uv slices}, which show $u(x,0.5)$ along the horizontal centerline and $v(0.5,y)$ along the vertical centerline. Along both centerlines, the predictions from the misspecified model show visible deviations from the reference solution, while those from the corrected model remain much closer to it.
	
	\begin{table}[h]
		\centering
		\begin{tabular}{c|c|c|c}
			\hline
			\hline
			& Misspecified model & Ensemble PINNs & Our approach \\
			\hline
			Error of $u$ & 0.0981 & 0.0131 & 0.0096 \\  
			\hline
			Error of $v$ & 0.1634 & 0.0223 & 0.0119\\
			\hline
			\hline
		\end{tabular}
		\caption{Non-Newtonian cavity flow. The errors are computed as the relative $L_2$ error of the predicted values for $u$ and $v$ compared to the reference solutions.}\label{tab:cavity_flow}
	\end{table}

	\section{Summary}\label{Summary}
	In this work, we developed a latent representation learning based model-correction framework by extending LVM-GP for PDE learning under misspecified physics with uncertainty quantification.
	The proposed method uses a confidence-aware encoder to construct a shared stochastic latent representation, and then employs a solution decoder together with a correction decoder to infer both the solution distribution and an associated discrepancy term.
	By conditioning on this shared latent representation, the two decoders jointly quantify uncertainty in both the solution and the correction under soft physics-informed constraints with noisy and sparse observations, while an auxiliary latent-space regularization controls the learned representation and improves robustness.
	Numerical experiments on an ODE system, a reaction-diffusion equation, and two-dimensional channel and cavity flows demonstrate that the proposed approach effectively mitigates model misspecification while providing uncertainty estimates for both the solution and the discrepancy.
	
	Future work will focus on improving the expressiveness and training efficiency of the correction network, as well as evaluating the proposed framework in broader application domains such as fluid dynamics, weather forecasting, and engineering simulations. In addition, extending LVM-GP to multiscale problems remains an important direction for future research.
	
	\FloatBarrier

	\bibliography{references}
	\newpage
	\appendix
	\setcounter{remark}{0}
	\setcounter{theorem}{0}
	\setcounter{definition}{0}
	\setcounter{figure}{0}
	\setcounter{table}{0}
	\renewcommand{\thetheorem}{\Alph{section}.\arabic{theorem}}
	\renewcommand{\thedefinition}{\Alph{section}.\arabic{definition}}
	\renewcommand{\theremark}{\Alph{section}.\arabic{remark}}
	\renewcommand{\thefigure}{\Alph{section}.\arabic{figure}}
	\renewcommand{\thetable}{\Alph{section}.\arabic{table}}
	
	\section{Additional results for ODE system}\label{Additional_results_for_ode_system}

\begin{figure}[H]
		\begin{center}
			\begin{overpic}[width=0.33\textwidth, trim=0 0 0 0, clip=True]{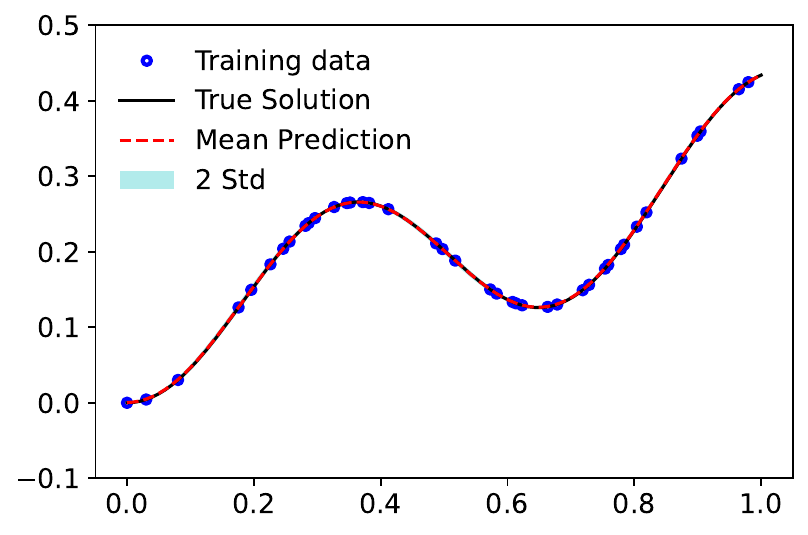}
				\put(52.5,67.5){{$u$}}
			\end{overpic}
			\begin{overpic}[width=0.33\textwidth, trim=0 0 0 0, clip=True]{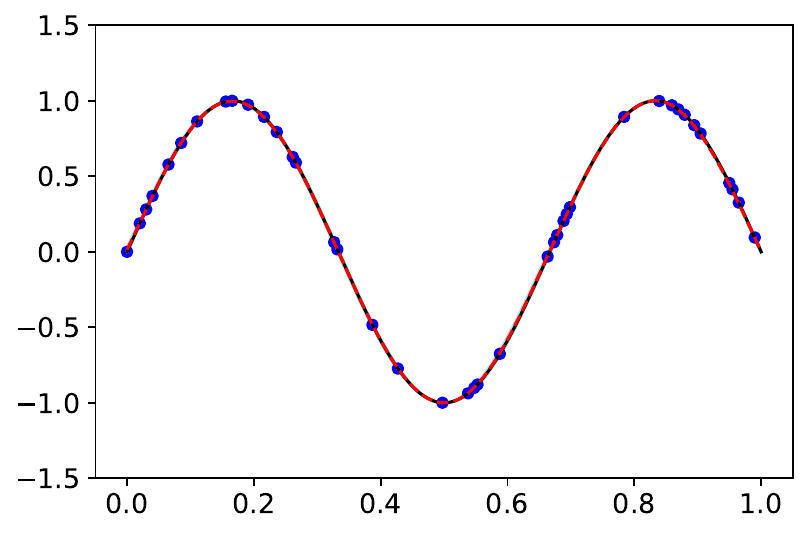}
				\put(52.5,67.5){{$f$}}
			\end{overpic}
			\begin{overpic}[width=0.33\textwidth, trim=0 0 0 0, clip=True]{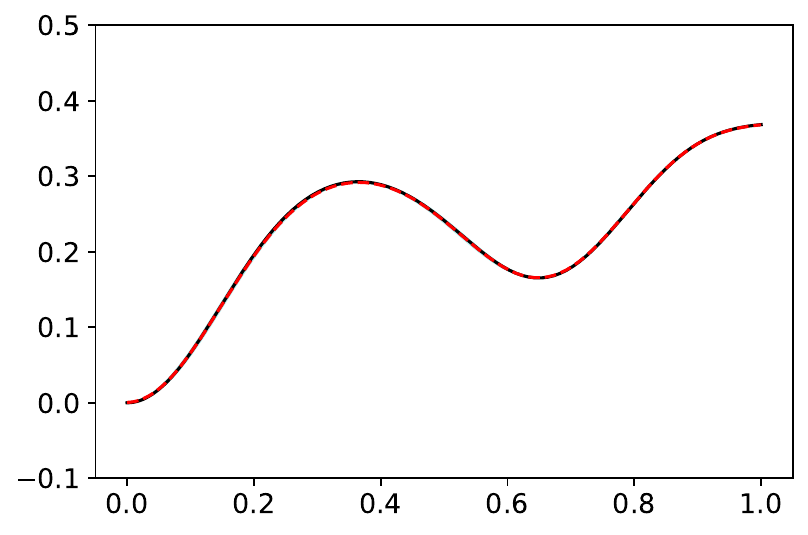}
				\put(50,67.5){{$\phi$}}
			\end{overpic}
			\subcaption{S1: LVM-GP with the correct physical model}
			
			\begin{overpic}[width=0.33\textwidth, trim=0 0 0 0, clip=True]{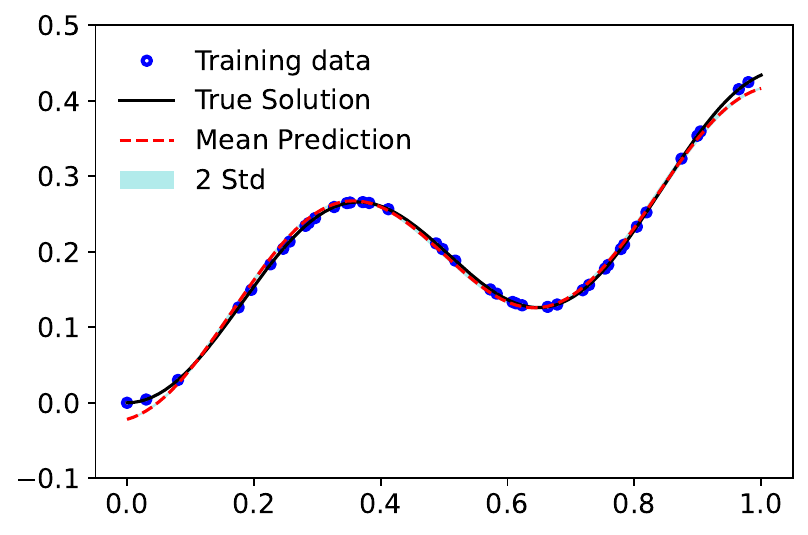}
			\end{overpic}
			\begin{overpic}[width=0.33\textwidth, trim=0 0 0 0, clip=True]{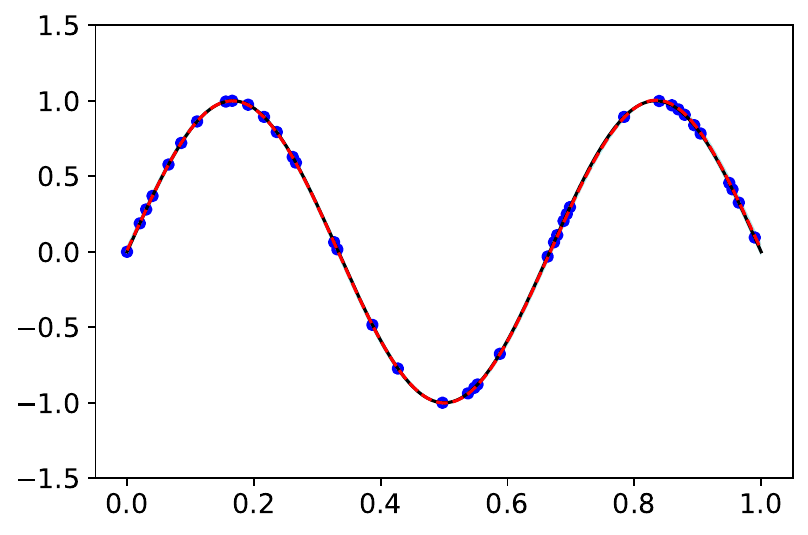}
			\end{overpic}
			\begin{overpic}[width=0.33\textwidth, trim=0 0 0 0, clip=True]{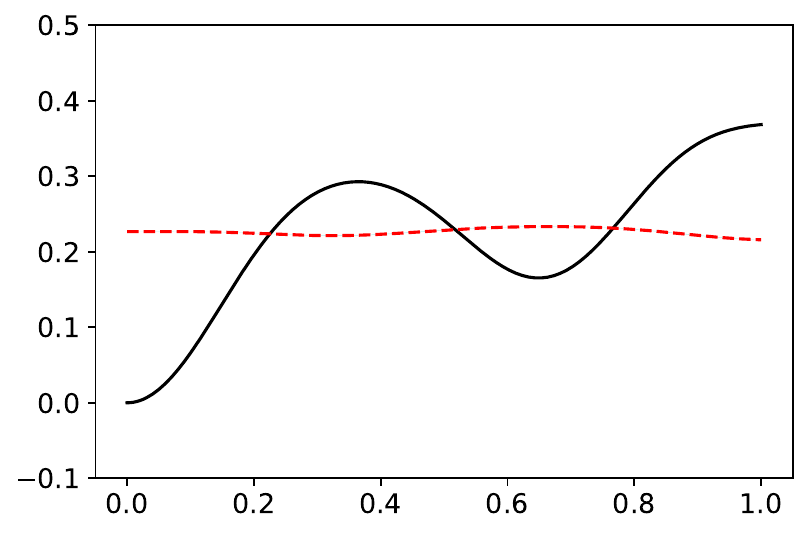}
			\end{overpic}
			\subcaption{S2: LVM-GP with a misspecified physical model}
			
			\begin{overpic}[width=0.33\textwidth, trim=0 0 0 0, clip=True]{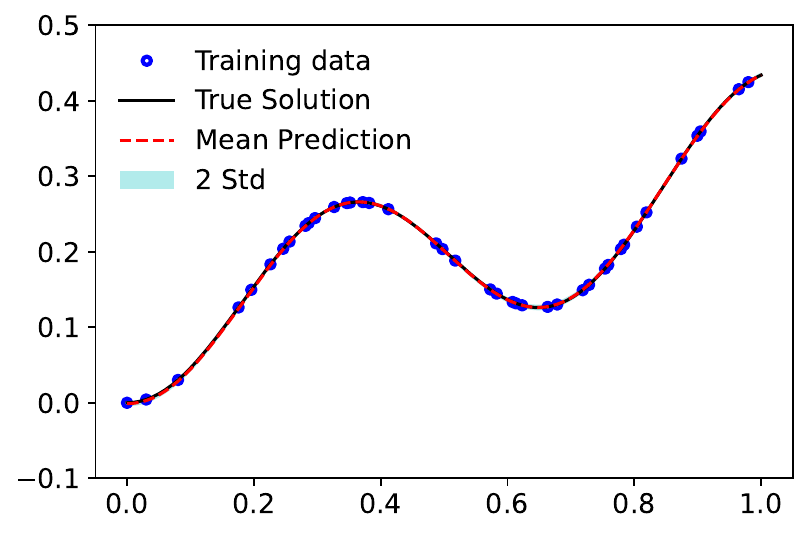}
			\end{overpic}
			\begin{overpic}[width=0.33\textwidth, trim=0 0 0 0, clip=True]{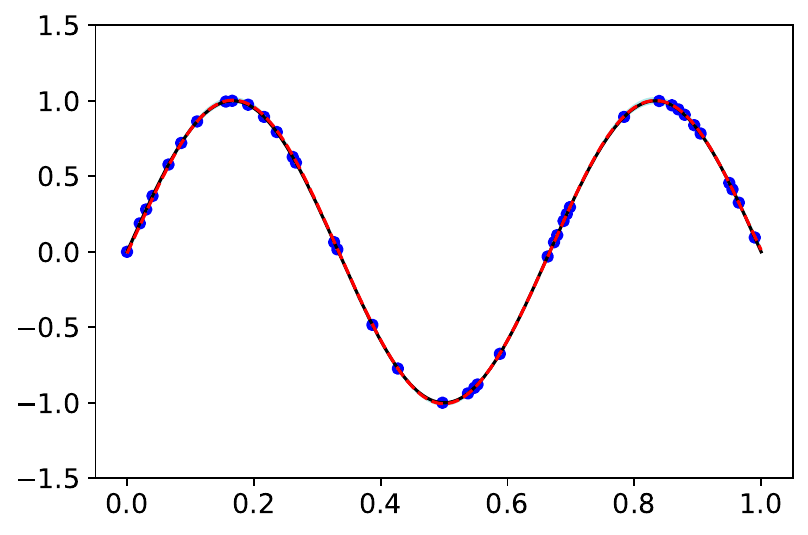}
			\end{overpic}
			\begin{overpic}[width=0.33\textwidth, trim=0 0 0 0, clip=True]{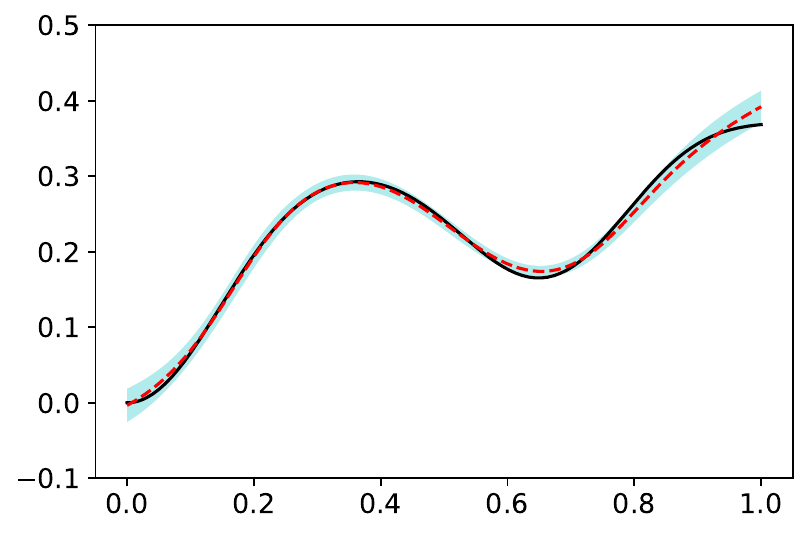}
			\end{overpic}
			\subcaption{S3: Correcting the misspecified physical model by a correction network using LVM-GP}
		\end{center}
		\caption{\textbf{ODE system}. Predictions of $u$, $f$, and $\phi$ under different scenarios using the same sparse and noise-free dataset. }
		\label{ode_noisefree_gappydata}
	\end{figure}

	As a supplement to Section~\ref{ODE_System}, we further investigate the case of sparse and noise-free data, where the time points $t$ are randomly sampled from the interval $t \in [0, 1]$. The results are shown in Figure~\ref{ode_noisefree_gappydata} and Table~\ref{tab:ode_gappy_noisefree}. The LVM-GP with correctly specified physical models achieves the best performance in both prediction accuracy and reliability (S1). The purely misspecified model yields the worst results (S2); however, after applying model correction (S3), the prediction accuracy is significantly improved, and all prediction errors fall within the uncertainty bounds, confirming the reliability of the corrected model. These results demonstrate the general applicability of our proposed method. Compared to the results in Section~\ref{ODE_System} where the data were dense and clean, the errors and uncertainties here are notably higher due to data sparsity and irregular sampling.
	
	\begin{table}[H]
		\centering
		\begin{tabular}{c|c|c|c|c|c}
			\hline
			\hline
			& $\tilde{\lambda}$ & Error of $\phi$ & Error of $u$ & Error of $f$ & Error of $\tilde{u}$ \\
			\hline
			S1: Known model & 1.4937 $\pm$ 0.0072 & 0.0029 & 0.0019 & 0.0050 & 0.0097 \\
			\hline
			S2: Misspecified model & 0.2321 $\pm$ 0.0036 & 0.4044 & 0.0303 & 0.0064 & 0.0990 \\
			\hline
			S3: Corrected model & 0.2 & 0.0278 & 0.0033 & 0.0081 & 0.0021 \\
			\hline
			\hline
		\end{tabular}
		\caption{\textbf{ODE system}. Predictions of $\phi, u, f$ and the reconstructed state $\tilde{u}$ are presented, with sparse and noise-free data. The table presents the relative $L_2$ errors between the reference solution and the mean predictions from LVM-GP (S1 and S2) as well as from the LVM-GP-based model correction (S3).}
		\label{tab:ode_gappy_noisefree}
	\end{table}
	
\end{document}